\documentclass[12pt]{amsart}
\advance\textheight\topskip
\usepackage[rightcaption]{sidecap}
\usepackage{amsmath,amsfonts,amsthm,amssymb,amscd}
\usepackage[mathscr]{eucal}

\usepackage[arrow]{xy}
\allowdisplaybreaks

\usepackage{hyperref}
\usepackage{times}

\numberwithin{equation}{section}
\makeatletter
\@addtoreset{equation}{section}

\def\@secnumfont{\bfseries}
\def\subsubsection{\@startsection{subsubsection}{3}%
  \z@{.5\linespacing\@plus.7\linespacing}{-.5em}%
  {\normalfont\bfseries}}
\def\paragraph{\@startsection{paragraph}{4}%
  \z@\z@{-\fontdimen2\font}%
  \normalfont\bfseries}
\def\subparagraph{\@startsection{subparagraph}{5}%
  \z@\z@{-\fontdimen2\font}%
  \normalfont\bfseries}

\makeatother

\newcommand{\xunderset}[2]{\underset{(#2)}{#1}}

\newlength{\iwidth}
\newcommand{\overwrite}[2]{\settowidth{\iwidth}{#2}%
  \smash[t]{\overset{\scriptscriptstyle#1}{#2}}\kern-.1em\kern-\iwidth\phantom{#2}}

\newcommand{\cI}{\mathcal{I}}

\newcommand{\ideal}{\mathscr{I}}

\newcommand{\Grass}[1]{\varepsilon(#1)}
\newcommand{\card}{\mathop{\mathrm{card}}}
\newcommand{\manX}{\mathscr{X}}

\newcommand{\func}[1]{C^\infty(#1)}
\newcommand{\sS}{\mathscr{S}}

\newcommand{\noo}[1]{{\boldsymbol{:}}\kern1pt #1{\boldsymbol{:}}}
\newcommand{\intpart}[1]{\left[#1\right]}
\newcommand{\dd}{\partial}

\newcommand{\bref}[1]{\textbf{\ref{#1}}}
\newcommand{\Lie}{\opfont{L}^{\phantom{y}}}
\newcommand{\ii}{\mathop{\opfont{I}}\nolimits^{\phantom{y}}}
\newcommand{\ghP}{\Bar{\mathscr{P}}}
\newcommand{\Hom}{\mathop{\mathrm{Hom}}\nolimits}
\newcommand{\contr}[2]%{\boldsymbol{\iota}^{\phantom{y}}_{#1}}
{\mbox{\kern1pt\large$\iota$}\!\left(#1\right)_{#2}}
\newcommand{\Contr}[2]{\left\langle#1,#2\right\rangle}

\newcommand{\Gtransf}{\mathscr{G}}

\newcommand{\Utransf}{\mathscr{U}}
\newcommand{\cH}{\mathcal{H}}
\newcommand{\Xterm}{\mathcal{X}}
\newcommand{\Yterm}{\mathcal{Y}}
\newcommand{\Zterm}{\mathcal{Z}}

\newcommand{\Sym}[2]{\mathop{\opfont{S}^{(#1)}}\limits_{#2}}

\newcommand{\ad}[2]{%
  {\smash[b]{\mathop{\opfont{ad}}\limits_{#1}}}{}^{\phantom{y}}_{#2}}

\newcommand{\gh}{\mathop{\mathrm{gh}}\nolimits}

\newcommand{\AQ}{\mathbb{A}^{\!\hbar}}
\newcommand{\AQgh}{\mathbb{A}^{\hbar}_{\mathrm{gh}}}

\newcommand{\PP}{\kern1pt\opfont{P}^{\phantom{y}}}

\newcommand{\thalf}{\tfrac{1}{2}}
\newcommand{\half}{\frac{1}{2}}
\newcommand{\uu}{\opfont{U}}

\newcommand{\vv}{\mathop{\opfont{V}}\nolimits}

\newcommand{\bull}[1]{\bullet{#1}\kern1pt}

\newcommand{\opfont}{\mathsf}
\newcommand{\algA}{\mathscr{A}}
\newcommand{\bA}{\Bar{\algA}}

\newlength{\Zlength}
\newcommand{\Z}[2][\relax]{%
  \settowidth{\Zlength}{\textsf{Z}}\smash[b]{%
    {\opfont{Z}\rule[-7pt]{0pt}{7pt}}^{#2}_{\kern-\Zlength #1}}}

\newcommand{\id}{\opfont{id}}
\newcommand{\idd}[1]{\id^{\tensor #1}}
\newcommand{\mult}{\opfont{m}}
\newcommand{\tensordots}{\tensor\dots\tensor}
\newcommand{\bb}{\opfont{b}}

\newcommand{\sss}{\mathop{\opfont{s}}\nolimits}

\newcommand{\oC}{\mathbb{C}}
\newcommand{\oZ}{\mathbb{Z}}
\newcommand{\tensor}{\mathbin{\otimes}}
\newcommand{\BRST}{\boldsymbol{\Omega}}

\newcommand{\ldot}{\mathbin{\boldsymbol{.}}}

\swapnumbers
\newtheorem{Thm}[subsection]{Theorem}
\newtheorem{thm}[subsubsection]{Theorem}

\newtheorem{lemma}[subsubsection]{Lemma}

\theoremstyle{definition}

\newtheorem{dfn}[subsubsection]{Definition}
\newtheorem{rem}[subsubsection]{Remark}

\begin{document}

%\begin{frontmatter}

\title%[Associativity in the BFV Formalism]%
[BRST on the Hochschild chain complex] {\vspace*{-4\baselineskip}
  \mbox{}\hfill\texttt{\small\lowercase{math}.QA/0301136}\\[3\baselineskip]
  Associativity and Operator Hamiltonian Quantization of Gauge
  Theories}

\author[I\,A\,Batalin]{I.\,A.~Batalin}
\author[A\,M\,Semikhatov]{A.\,M.~Semikhatov}
\address{\mbox{}\kern-\parindent Tamm Theory Division, Lebedev Physics
  Institute, Russian Academy of Sciences}

\begin{abstract}
  We show that the associative algebra structure can be incorporated
  in the BRST quantization formalism for gauge theories such that
  extension from the corresponding Lie algebra to the associative
  algebra is achieved using operator quantization of reducible gauge
  theories.  The BRST differential that encodes the associativity of
  the algebra multiplication is constructed as a second-order
  quadratic differential operator on the bar resolution.
\end{abstract}

\maketitle
\thispagestyle{empty}
%% \flushcolumns

\renewcommand{\contentsname}{}
\setcounter{tocdepth}{1}
\vspace*{-36pt}

{\addtolength{\parskip}{-4pt}
  \scriptsize
  \tableofcontents}
\enlargethispage{18pt}

% \begin{center}
%  \parbox{\textwidth}{
%    \begin{multicols}{2}
%      {\footnotesize
%        \tableofcontents}
%    \end{multicols}
%    }
% \end{center}

\section{Introduction}
\subsection{}  BRST quantization methods for constrained systems
(synonym\-ous\-ly, gauge theories)~\cite{[BV],[BF-Ham],[HT]} are
recognized as a powerful approach reaching beyond the contexts in
which it was originally created.  In the Hamiltonian formalism, the
BRST quantization of a first-class constrained system amounts to
constructing an odd differential (the ``BRST'' operator) in the
constrained system extended by ghosts and the conjugate momenta (which
are recognized as Koszul--Tate variables at the classical
level~\cite{[S-gh]}).  But the formalism has been limited to algebraic
structures built on (graded) \textit{antisymmetric} operations
(generated by commutators or Poisson brackets).

In this paper, we show that the Hamiltonian BRST formalism is also
applicable to graded \textit{associative} algebras and is therefore
not limited to Lie-like structures.  For a given associative
algebra~$\algA$, we construct a differential $\BRST$ such that the
relation $\BRST^2=0$ is a ``BRST encoding'' of the associativity of
multiplication in~$\algA$.

The construction has a ``trivial'' part, the associated (graded) Lie
algebra $(\algA,[\;,\,])$, which is treated by the standard BRST
methods, and a ``difficult'' part, the extension of the BRST scheme to
the associative algebra.  This extension is constructed using the
machinery of \textit{reducible} gauge theories~\cite{[BF-red]}, with
the bar resolution of~$\algA$ viewed as the data defining a reducible
gauge theory.  In accordance with the BRST ideology, we introduce
ghosts for each term in the bar resolution, quantize them, and seek
the BRST differential~$\BRST$ of ghost number~$1$; we also require
$\BRST$ to be at most quadratic in the ghosts and at most bilinear in
the momenta, with no individual momentum being squared in~$\BRST$.
The BRST differential is therefore a quadratic second-order
\textit{differential operator on the bar resolution}; moreover, its
purely quadratic part is an operator with scalar coefficients.  The
exact statement is formulated in~\bref{thm:Omega2} below.

\subsection{} Although our construction of the differential~$\BRST$ is
motivated by BRST methods, similarities with the known BRST formalism
for reducible gauge theories are somewhat limited because the bar
resolution is infinite, and the corresponding reducible gauge theory
is therefore of an infinite reducibility degree.  Despite some
effort~\cite{[BV-red],[BF-red],[Vilk]}, additionally motivated by
possible string theory applications~\cite{[KTvP]}, infinitely
reducible gauge theories have not been given a complete formulation
that would allow proving the existence of the BRST differential (or a
solution to the master equation in the Lagrangian formulation).  The
standard inductive argument applicable to finitely reducible gauge
theories fails in the infinitely reducible case, and we must therefore
independently prove the existence in the associative algebra setting;
for this, we construct a recursive procedure that yields a solution of
the equation~$\BRST^2=0$.  This also bypasses another complication: in
contrast to ``genuine'' gauge theories, where first-class constraints
satisfy the so-called involution relations (which become
Poisson-bracket relations as~$\hbar\to0$)
\begin{equation}\label{eq:involution}
  [T_\alpha,\,T_\beta]
  = i\hbar\sum_{\gamma} U_{\alpha\beta}^\gamma T_\gamma
\end{equation}
with some operators~$U_{\alpha\beta}^\gamma $, the commutators
in~$\algA$ do not involve a Planck constant~$\hbar$.  In other words,
the Planck constant is equal to~$1$ (in fact, to~$-i$) in the
formalism proposed here.  Hence, there is no classical limit of the
corresponding gauge theory, and the BRST formalism must be applied
directly at the operator level.  Another consequence of $\hbar$ being
equal to~$1$ is a problem (which we entirely ignore) in interpreting
infinite series that are no longer formal series in~$\hbar$.

\subsection{} The paper is organized as follows.  For completeness, 
we recall basic facts about reducible gauge theories in
Sec.~\ref{sec:reducible}.
%% (most part of it may be skipped by some of the readers, in particular
%% because the known gauge-theory facts do not directly imply our main
%% result, the existence of the BRST differential).  
In~\bref{sec:bar}, we explain the relation between
Sec.~\ref{sec:reducible} and the main part of Sec.~\ref{sec:Omega}.
We introduce ghosts and construct a differential in~\bref{sec:BRST} as
the BRST operator in a specific (infinitely) reducible gauge theory.
The main result that $\BRST^2=0$ is formulated in~\bref{thm:Omega2}
and is proved in~\bref{sec:proof}--\bref{sec:end-proof-thm}.  Several
additional remarks are collected in Sec.~\ref{sec:other}.

\section{Reducible Constrained Systems} \label{sec:reducible}
In this section, we summarize the main points of the Hamiltonian BRST
quantization of constrained systems.  We use it as a motivation for
what follows, even though direct application of the basic BRST
formalism to the associative algebra case is hindered by the infinite
reducibility and the absence of the classical limit.  The actual link
between the BRST formalism and our main construction is explained
in~\bref{sec:bar}.

\subsection{Classical reducible constrained systems}
\addcontentsline{toc}{section}{\numberline{\thesubsection.}Classical
  reducible constrained systems} 
A classical \textit{reducible constrained system}~\cite{[BF-red]}
consists of
\begin{itemize}
\item a symplectic manifold $\manX$ (in fact, a symplectic vector
  space) and (sufficiently smooth) functions $T_{\alpha_0}$,
  $\alpha_0\in\cI_0$, on~$\manX$, called \textit{constraints}, whose
  zero locus is called the constraint ``surface''~$\sS\subset\manX$
  and is viewed by physicists as something very close to a smooth
  manifold,
  
\item the functions $Z_{\alpha_n}^{\alpha_{n-1}}$, $\alpha_i\in\cI_i$,
  on~$\manX$ satisfying the rank assumption and ``zero-mode''
  equations~\eqref{ZT} and~\eqref{ZZ} given below.

\end{itemize}

We now consider these ingredients in more detail.  If the given
constraints are linearly independent over $\func{\manX}$, the theory
is said to be irreducible; if they are not, there exist functions
$\overwrite{1}{Z}_{\alpha_1}^{\alpha_0}$, $\alpha_1\in\cI_1$,
on~$\manX$ such that
\begin{equation}\label{ZT}
  \sum_{\alpha_0\in\cI_0}
  \overwrite{1}{Z}_{\alpha_1}^{\alpha_0} T_{\alpha_0}=0.
\end{equation}
The functions $\overwrite{1}{Z}_{\alpha_1}^{\alpha_0}$ can in turn be
linearly dependent over $\func{\manX}$, which gives rise to
$\overwrite{2}{Z}_{\alpha_2}^{\alpha_1}$, and so on.  More precisely,
each subsequent linear dependence is only required in a ``weak'' form,
i.e., modulo the ideal $\ideal$ generated by~$\{T_{\alpha_0}\}$,
\begin{equation}\label{ZZ}
  \sum_{\alpha_{n-1}}  
  \overwrite{n}{Z}_{\alpha_n}^{\alpha_{n-1}}
  \overwrite{n-1}{Z}_{\alpha_{n-1}}^{\alpha_{n-2}}
  \in\ideal,
  \qquad \alpha_i\in\cI_i,\quad n\geq2
\end{equation}
(we often omit top labels $n$ from the notation in what follows).
Upon restriction to the constraint surface, this gives a complex,
which is exact by definition.\footnote{The entire complex, including a
  specific choice of~$\overwrite{n}{Z}$, makes part of the definition
  of a reducible constrained system (reducible gauge theory): although
  the exact sequence furnished by $\overwrite{n}{Z}$ splits in
  physical applications, a given splitting is not canonical and in
  realistic theories, moreover, typically violates some important
  symmetries (e.g., the Lorentz covariance) or locality.}  With each
$\overwrite{n}{Z}_{\alpha_n}^{\alpha_{n-1}}$ viewed as a rectangular
matrix, the ranks of their restrictions to~$\sS$ must therefore
satisfy
\begin{align*}
  \mathop{\mathrm{rank}}
  \bigl(\overwrite{n+1}{Z}_{\alpha_{n+1}}^{\alpha_{n}}\bigr)
  \Bigr|_{\sS}
  + \mathop{\mathrm{rank}}
  \bigl(\overwrite{n}{Z}_{\alpha_n}^{\alpha_{n-1}}\bigr)
  \Bigr|_{\sS}
  ={}& \card\cI_n,\quad n\geq1,\\
  \intertext{and in addition,}
  \mathop{\mathrm{rank}}
  \bigl(\overwrite{1}{Z}_{\alpha_1}^{\alpha_{0}}\bigr)
  \Bigr|_{\sS} +
  \mathop{\mathrm{rank}}\bigl(\tfrac{\dd T_{\alpha_0}}{\dd x^i}
  \bigr)\Bigr|_{\sS}
  ={}& \card\cI_0,
\end{align*}
where $\{x^i\}$ is any local coordinate system in a neighborhood of
$\sS$ in~$\manX$.  \textit{These ranks are assumed to be constant in
  some neighborhood of~$\sS$ in~$\manX$~\cite{[BF-red],[HT]}}.  We
generally refer to $\overwrite{n}{Z}$ as ``zero modes.''

\begin{dfn}
  A constrained system is said to be $\ell$-reducible if
  $\overwrite{\ell+1}{Z}=0$, but $\overwrite{\ell}{Z}\not\equiv0$.
  In particular, a $0$-reducible theory is irreducible
  (relations~\eqref{ZT} are already absent).
\end{dfn}

%% \subsubsection{} Different terms of the complex can have arbitrary
%% $\oZ_2$ gradings, and the gradings of $Z_{\alpha_{n}}^{\alpha_{n-1}}$
%% are determined appropriately.  We do not introduce special notation
%% for the linear spaces (in fact, bundles over~$\sS$) between which
%% $\overset{n}{Z}$ act, and specify only the $\oZ_2$ gradings of
%% $Z_{\alpha_{n}}^{\alpha_{n-1}}$ as
%% \begin{equation}\label{Z2}
%%   \Grass{Z_{\alpha_{n}}^{\alpha_{n-1}}} = \varepsilon_{\alpha_n}
%%   + \varepsilon_{\alpha_{n-1}}\quad\text{and}\quad
%%   \Grass{T_{\alpha_0}} = \varepsilon_{\alpha_0}.
%% \end{equation}

\subsubsection{} A constrained system can be extended by auxiliary
variables, called ghosts and the conjugate momenta, such that there
exists an odd Hamiltonian vector field $\{\Omega,{-}\}$ (the BRST
operator) whose square is zero and the lowest terms in the expansion
of $\Omega$ in the ghosts involve the constraints and the zero modes
$\overwrite{n}{Z}$.  To avoid repetition, we consider the details in
the quantum setting.

\subsection{Quantum reducible constrained systems}
\addcontentsline{toc}{section}{\numberline{\thesubsection.}{Quantum
    reducible constrained systems}} The mathematically rigorous
existence of the quantum theory is a subtle issue, and the following
statements may be sensitive to the chosen quantization.\pagebreak[3]
Deformation quantization~\cite{[BFFLS],[F],[K]} alone does not
automatically allow speaking of \textit{operator} relations, but we
proceed in terms of these to recapitulate the basic gauge-theory
folklore.

\subsubsection{} Quantum-mechanically, $T_{\alpha_0}$ and
$Z_{\alpha_n}^{\alpha_{n-1}}$ become operators, i.e., elements of an
algebra~$\AQ$.  Equation~\eqref{ZT} then retains its form in terms of
elements of~$\AQ$, and Eqs.~\eqref{ZZ} become
\begin{equation}\label{ZZ-q}
  \sum_{\alpha_{n-1}}  
  \overwrite{n}{Z}_{\alpha_n}^{\alpha_{n-1}}
  \overwrite{n-1}{Z}_{\alpha_{n-1}}^{\alpha_{n-2}}
  =
  \sum_{\alpha_0}
  \Pi_{\alpha_n}^{\alpha_{n-2}\alpha_0} T_{\alpha_0} + \hbar \AQ,
  \qquad \alpha_i\in\cI_i,
\end{equation}
with some $\Pi_{\alpha_n}^{\alpha_{n-2}\alpha_0}\in\AQ$ (the relations
$\overwrite{n}{Z}\cdot\overwrite{n-1}{Z}\equiv0\;\mathrm{mod}\;\ideal$
are therefore reproduced only as $\hbar\to0$).

\subsubsection{Ghost content} \label{sec:ghosts} The algebra $\AQ$ is
extended to $\AQgh$ by a set of operators $\{C^A\}$ (ghosts) and
$\{\ghP_A\}$ (conjugate momenta) satisfying the canonical graded
commutation relations
\begin{equation*}
  [\ghP_B, C^A]\equiv \ghP_B C^A
  - (-1)^{\Grass{\ghP_B} \Grass{C^A}}C^A \ghP_B
  = i\hbar\,\delta_B^A.
\end{equation*}
Here, 
\begin{itemize}
\item $A$ is a collection of (multi)indices,
  $A=\{\alpha_0,\alpha_1,\dots,\alpha_{\ell}\}$, $\alpha_i\in\cI_i$,
  and the ghosts are therefore a collection $\{C^A\} =
  \{C^{\alpha_0},\dots,C^{\alpha_\ell}\}$, $\alpha_i\in\cI_i$, and
  similarly for the momenta, $\{\ghP_B\} =
  \{\ghP_{\beta_0},\dots,\ghP_{\beta_\ell}\}$, $\beta_i\in\cI_i$.
  
\item $\oZ_2$-gradings of the ghosts and the momenta are
  \begin{equation*}
    \Grass{C^{\alpha_n}}=\Grass{\ghP_{\alpha_n}}
    =\varepsilon_{\alpha_n}+n+1;
  \end{equation*}
  
\item ghost-number assignments for the ghosts and the momenta are
  \begin{equation*}
    \gh C^{\alpha_n}=n+1,\qquad
    \gh \ghP_{\alpha_n} = -n-1.
  \end{equation*}
\end{itemize}
  
The following statement gives an operational tool for the Hamiltonian
quantization of constrained systems (see comments
in~\bref{rem:existence} in what follows).
\begin{thm}[Batalin--Fradkin] In a quantum $\ell$-reducible constrained
  system, there exists a $2$-nilpotent odd operator~$\Omega$ with
  ghost number~$1$ (the BRST differential),
  \begin{equation*}
    \Grass{\Omega}=1,\qquad \gh\Omega=1,\qquad \Omega^2=0,
  \end{equation*}  
  of the form
  \begin{equation}\label{Omega-gen}
    \Omega=\sum_{n\geq1}\sum_{m\geq0}
    \sum_{\substack{A_1,\dots,A_m\\
        B_1,\dots,B_n}}
    C^{A_m}\dots C^{A_1} V_{A_1\dots A_m}^{B_n\dots B_1}
    \ghP_{B_1}\dots\ghP_{B_n}\in\AQgh,
  \end{equation}  
  where $V_{A_1\dots A_m}^{B_n\dots B_1}$ are operators in~$\AQ$ such
  that
  \begin{equation}\label{boundary}
    V_{\alpha_0}=T_{\alpha_0},\qquad
    V_{\alpha_n}^{\alpha_{n-1}}
    =\overwrite{n}{Z}_{\alpha_n}^{\alpha_{n-1}},
  \end{equation}
  and as $\hbar\to0$, \ $T_{\alpha_0}$ and
  $\overwrite{n}{Z}_{\alpha_n}^{\alpha_{n-1}}$ become the respective
  data specifying a classical $\ell$-reducible gauge theory and
  $\Omega$ becomes an odd function satisfying the Poisson-bracket
  relation $\{\Omega,\Omega\}=0$ (with the Poisson bracket between the
  ghosts and the momenta determined by the above commutator).
\end{thm}

\subsubsection*{Remarks}

\subsubsection{} \label{rem:existence} 
The theorem claims the existence of the operators $V_{A_1\dots
  A_m}^{B_n\dots B_1}$ in addition to those fixed by ``boundary
conditions''~\eqref{boundary}.  In other words, given the constraints
and the $\overwrite{i}{Z}$ operators, a nilpotent $\Omega$ with the
ghost number~$1$ can be constructed as a formal series in the ghosts
and momenta starting with the lower-order terms
\begin{gather}\label{bc}
  \Omega = \sum_{\alpha_0}C^{\alpha_0}T_{\alpha_0}
  + \sum_{n=1}^{\ell}\sum_{\alpha_n,\alpha_{n-1}}
  C^{\alpha_n}
  \overwrite{n}{Z}_{\alpha_n}^{\alpha_{n-1}}\ghP_{\alpha_{n-1}}
  +\dots.
\end{gather}

The existence of a quantum BRST operator given by a series in the
ghosts has never been proved explicitly, however.  Its classical
counterpart is known to exist if $\manX$ is a symplectic linear space
and the homology of the Koszul--Tate operator is concentrated in ghost
number zero (this is the actual significance of the rank assumptions).
The classical BRST operator can then be subjected to deformation
quantization, and its quantization can be shown to exist in terms of
symbols (formal series in~$\hbar$) if the classical operator
$\{\Omega,{-}\}$ has trivial homology in ghost number~$1$.

\subsubsection{} 
%% In accordance with  $\oZ_2$-gradings of the ghosts,
The coefficients in~\eqref{Omega-gen} 
%% have certain symmetry properties: writing
are usually redefined as
\begin{equation*}
  V_{A_1\dots A_m}^{B_n\dots B_1}=
  (-1)^{E_{A_1\dots A_m}^{B_n\dots B_1}}\,
  \tfrac{1}{m!n!}\,U_{A_1\dots A_m}^{B_n\dots B_1},
\end{equation*}
%% (with conventional factorials), 
where the sign factors defined by
\begin{multline*}
  E_{A_1\dots A_m}^{B_n\dots B_1}=
  \sum_{k=1}^{\intpart{\frac{m}{2}}}\Grass{C^{A_{2k}}}
  + \sum_{k=1}^{\intpart{\frac{n}{2}}}\Grass{\ghP_{B_{2k}}}\\*
  + \intpart{\tfrac{m}{2}}\max\gh(C^{A_1},\dots,C^{A_m})
  + \intpart{\tfrac{n}{2}}\max\gh(\ghP_{B_1},\dots,\ghP_{B_n})
\end{multline*}
(with $\intpart{~}$ denoting the integer part) are chosen such that
\begin{align*}
  U_{A_1\dots A_kA_{k+1}\dots A_m}^{B_n\dots B_1}={}&
  -(-1)^{(\varepsilon_{A_k}+1)(\varepsilon_{A_{k+1}}+1)}
  U_{A_1\dots A_{k+1}A_{k}\dots A_m}^{B_n\dots B_1},\\
  U_{A_1\dots A_m}^{B_n\dots B_{k+1}B_k\dots B_1}={}&
  -(-1)^{(\varepsilon_{B_k}+1)(\varepsilon_{B_{k+1}}+1)}
  U_{A_1\dots A_m}^{B_n\dots B_{k}B_{k+1}\dots  B_1},
\end{align*}
where $\varepsilon_{A_k}=\Grass{C^{A_k}}$.  
%% The signs $(-1)^{E_{A_1\dots A_m}^{B_n\dots B_1}}$ are conventional
%% and are chosen such that 
This minimizes the number of explicit sign factors when the relations
generated by Eqs.~\eqref{eq:involution} and the reducibility relations
are written in terms of $U_{A_1\dots A_kA_{k+1}\dots A_m}^{B_n\dots
  B_1}$; (inevitable) sign factors then occur only in
cyclic-permutation sums.

\subsubsection{} \label{sec:eq-list} In what follows, we specialize to
the case where $\varepsilon_{\alpha_i}=0$, and therefore the
$\oZ_2$-gradings of the ghosts are
\begin{equation*}
  \Grass{C^{\alpha_i}}=i+1,\qquad i=0,\dots,\ell,
\end{equation*}
and all $U_{A_1\dots A_m}^{B_n\dots B_1}$ are even,
$\Grass{U_{A_1\dots A_m}^{B_n\dots B_1}}=0$.  The lower equations
following from $\Omega\Omega=0$ are given by involution
relations~\eqref{eq:involution}, ``zero-mode'' conditions~\eqref{ZT},
and (with summation over repeated indices understood)
\begin{align}\label{eq:THREE}
  (i\hbar)^{-1}\,[T_{\alpha_0},Z_{\alpha_1}^{\beta_0}]
  - Z_{\alpha_1}^{\gamma_0}\,U_{\gamma_0\alpha_0}^{\beta_0}
  &{}=-U_{\alpha_0\alpha_1}^{\gamma_0\delta_0}
  \Pi_{\gamma_0\delta_0}^{\beta_0}
  - U_{\alpha_0\alpha_1}^{\beta_1}\,Z_{\beta_1}^{\beta_0},\\
  &\phantom{{}={}}
  2\Pi_{\alpha_0\beta_0}^{\gamma_0}
  =T_{[\alpha_0}\delta_{\beta_0]}^{\gamma_0}
  -i\hbar\,U_{\alpha_0\beta_0}^{\gamma_0},\notag\\
  Z_{\alpha_2}^{\alpha_1}\,Z_{\alpha_1}^{\alpha_0}&=
  -U_{\alpha_2}^{\beta_0\gamma_0}\,\Pi_{\beta_0\gamma_0}^{\alpha_0},\\ 
  \label{eq:SEVEN}
  (i\hbar)^{-1}\,
  \bigl[Z_{\alpha_1}^{[\alpha_0},\,Z_{\beta_1}^{\beta_0]}\bigr]
  - Z_{\{\alpha_1}^{\gamma_0}\,U_{\gamma_0\beta_1\}}^{\alpha_0\beta_0}
  &{}=-U_{\alpha_1\beta_1}^{[\alpha_0 
    \gamma_0\delta_0}\,
  \Pi_{\gamma_0\delta_0}^{\beta_0]}
  - U_{\alpha_1\beta_1}^{[\alpha_0\gamma_1}\,Z_{\gamma_1}^{\beta_0]}\\*
  &\phantom{{}={}}
  +U_{\alpha_1\beta_1}^{\alpha_2}\,U_{\alpha_2}^{\alpha_0\beta_0}
  +i\hbar\,U_{\alpha_1\beta_1}^{\gamma_0\gamma_1}\,
  U_{\gamma_0\gamma_1}^{\alpha_0\beta_0},\notag\\
  \label{eq:SIX}
  Z_{\{\alpha_1}^{\alpha_0}\,U_{\alpha_0\beta_1\}}^{\gamma_1}
  + U_{\alpha_1\beta_1}^{\alpha_2}\,Z_{\alpha_2}^{\gamma_1}
  &{}=-U_{\alpha_1\beta_1}^{\alpha_0\delta_1}
  \left(T_{\alpha_0}\,\delta_{\delta_1}^{\gamma_1}
    +i\hbar\,U_{\alpha_0\delta_1}^{\gamma_1}\right)
\end{align}
(plus infinitely many more relations).  Here, $A^{[\alpha\beta]} =
A^{\alpha\beta} - A^{\beta\alpha}$ and $A^{\{\alpha\beta\}} =
A^{\alpha\beta} + A^{\beta\alpha}$; whenever a group of indices are
bracketed, the antisymmetrization and symmetrization operations apply
to only the leftmost and rightmost indices.

\subsection{Reducible closed-algebra theories}
\label{sec:closed}
\addcontentsline{toc}{section}{\numberline{\thesubsection.}Reducible closed-algebra theories}
\begin{dfn}\label{def:closed}
  First-class constraints $T_{\alpha_0}$, $\alpha_0\in\cI_0$, are said
  to generate a \textit{closed gauge algebra} if $\Omega$ is at most
  quadratic in the ghosts.
\end{dfn}

The issue of a ``closed algebra'' is more subtle in the Hamiltonian
than in Lagrangian BRST formalism.  In what follows, we consider
closed-algebra theories satisfying an additional assumption that
\textit{the BRST differential $\Omega$ is at most bilinear in the
  momenta, with none individual momentum $\ghP_{\alpha_i}$ entering
  squared} (in reducible theories, this is not a consequence of the
closed-algebra condition, but in all known examples where both
conditions are satisfied, the corresponding gauge generators in the
Lagrangian formulation generate a closed gauge algebra).

\subsubsection{} 
The ``gauge theory'' considered in Sec.~\ref{sec:Omega} has a closed
gauge algebra but is infinitely reducible, and we therefore allow sums
over infinitely many ghosts.  We use indices $A$, $B$, \dots to label
\textit{all} ghosts (i.e., $C^{\alpha_0}$, $C^{\alpha_1},\dots$ in
terms of~\bref{sec:ghosts}).  In a closed-algebra theory, the BRST
differential can be written as
\begin{equation}\label{gen-reducible}
  \Omega=\sum_A C^A \tau_A
  + \thalf\sum_{A,B}(-1)^{\Grass{C^B}+1}C^BC^A U_{AB},
\end{equation}
where $\tau_A$ and $U_{AB}$ are \textit{functions of the ghost
  momenta~$\ghP_C$}.  In the notation of~\bref{sec:ghosts},
$\tau_{\alpha_0}=T_{\alpha_0}$ are obviously the original
constraints;\pagebreak[3] all $\tau_{\alpha_n}$ with $n\geq1$ are
$\ghP_C$-dependent. The equation $\Omega^2=0$ then amounts to
independent vanishing conditions for the terms of the first, second,
and third degree in $C^A$.  Anticipating the structure of the BRST
differential in Sec.~\ref{sec:Omega} (where $U_{AB}$ are
\textit{scalar-valued}), we write these equations in the case where
\begin{equation}\label{scalar-case}
  [U_{AB},U_{CD}]=0\quad\text{and}\quad [\tau_A,U_{BC}]=0.
\end{equation}
In the linear order in $C^A$, the equations are given by (with
summations over repeated indices understood)
\begin{gather} \label{C1}
  [\tau_A,C^B]\tau_B+\thalf[[\tau_A,C^C],C^B]
  U_{BC}(-1)^{\Grass{C^A}+1}=0.
\end{gather}
In the quadratic and cubic orders, the respective equations are
\begin{multline}\label{C2}
  [\tau_A,\tau_B]=[U_{AB},C^C]\tau_C - 
  [\tau_A,C^C]U_{CB}(-1)^{\Grass{C^B}+1}\\*
  +[\tau_B,C^C]U_{CA}(-1)^{(\Grass{C^A}+1)\Grass{C^B}}
  +\thalf[[U_{AB},C^D],C^C]U_{CD}(-1)^{\Grass{C^D}+1}
\end{multline}
and
\begin{equation} \label{C3}
  [U_{AB},C^D]U_{DC}(-1)^{(\Grass{C^A}+1)(\Grass{C^C}+1)} +
  \text{cycle}(A,B,C)=0.
\end{equation}

\section{The BRST Differential} \label{sec:Omega}
In this section, we construct the BRST differential~$\BRST$ as a
differential operator on the bar resolution of a given associative
algebra~$\algA$.  In~\bref{sec:BRST}, we introduce the necessary
structures in the physical context, in the guise of ghosts in a
reducible gauge theory.  The main result that $\BRST$ \textit{is} a
differential is formulated in Theorem~\bref{thm:Omega2}; a recursive
solution to $\BRST\,\BRST=0$ is constructed in~\bref{sec:proof} and
the proof that $\BRST\,\BRST=0$ extends to~\bref{sec:end-proof-thm}.

\addtocounter{subsection}{-1}
\subsection{Tensor algebra preliminaries} \label{sec:prelim}
Let $\algA$ be an associative graded algebra with a unit.  We write
\begin{align*}
  \mult_{1,i}:\algA^{\tensor n}&{}\to\algA^{\tensor(n-1)}\\*
  a_1\tensordots a_{n}&{}\mapsto
  (-1)^{\Grass{a_i}(\Grass{a_2}+\dots+\Grass{a_{i-1}})}
  a_1 a_{i}\tensor a_2\tensor\ldots\not{\kern-2pt a_i}
  \ldots\tensor a_{n}
\end{align*}
for $2\leq i\leq n$ and use the somewhat redundant notation
$\mult_{i,1}^*$ (with both the transposed subscripts and the asterisk)
for the reversed-order multiplication (as before, placed in the first
tensor factor)
\begin{align*}
  \mult_{i,1}^*:\algA^{\tensor n}&{}\to\algA^{\tensor(n-1)}\\*
  a_1\tensordots a_{n}&{}\mapsto
  (-1)^{\Grass{a_i}(\Grass{a_1}+\dots+\Grass{a_{i-1}})}
  a_{i}a_1\tensor a_2\tensor\ldots\not{\kern-2pt a_i}
  \ldots\tensor a_{n}.
\end{align*}
We also use the notation
\begin{align*}
  \ad{1}{a'}\ldot a_0\tensordots a_n &= [a',a_0]\tensor
  a_1\tensordots a_n,
\end{align*}
and let $\sss$ be the ``right shift'' mapping
\begin{align*}
  \sss:\algA^{\tensor n}&{}\to\algA^{\tensor(n+1)}\\*
  a_1\tensordots a_n &{}\mapsto 1\tensor a_1\tensordots a_n.
\end{align*}
Finally, we need the operation of moving the $i$th tensor factor to
the first position,
\begin{equation*}
  \PP_{1,i}: a_1\tensordots a_n
  \mapsto
  (-1)^{\Grass{a_i}(\Grass{a_1}+\dots+\Grass{a_{i-1}})}
  a_{i}\tensor a_1\tensor a_2\tensor\ldots\not{\kern-2pt a_{i}}
   \ldots\tensor a_n,
\end{equation*}
and the inverse operation
\begin{equation*}
  \PP_{i,1}: a_1\tensordots a_n
  \mapsto
  (-1)^{\Grass{a_1}(\Grass{a_2}+\dots+\Grass{a_{i}})}
  a_2\tensordots a_i\tensor a_1\tensor a_{i+1}\tensordots a_n.
\end{equation*}

\subsection{The bar resolution and the ``zero modes''} \label{sec:bar}
\addcontentsline{toc}{section}{\numberline{\thesubsection.}The bar
  resolution and the ``zero modes''}

We consider the bar resolution of the algebra~$\algA$ by free
$\algA$-bimodules
\begin{equation}\label{resolution}
  {}\dots\xrightarrow{\bb'}\algA\tensor\bA^{\tensor2}\tensor\algA
  \xrightarrow{\bb'}\algA\tensor\bA\tensor\algA
  \xrightarrow{\bb'}\algA\tensor\algA
  \xrightarrow{\bb'}\algA\to 0,
\end{equation}
where $\bA=\algA/\oC$ and the differential is given by
\begin{equation*}
  \bb' a_0\tensor a_1\tensordots a_n=
  \sum_{i=0}^{n-1}(-1)^i a_0\tensordots a_{i-1}
  \tensor a_i a_{i+1}\tensor a_{i+2}\tensordots a_n
\end{equation*}
(in particular, $\algA\tensor\algA\xrightarrow{\bb'}\algA$ is simply
the multiplication $\algA\tensor\algA\xrightarrow{\mult}\algA$).  This
complex is contractible with the contracting homotopy given by~$\sss$,
\begin{equation*}
  \sss\bb' + \bb'\sss = \id.  
\end{equation*}

To establish connection with Sec.~\ref{sec:reducible}, we define
linear mappings $\Z{n} : \algA^{\tensor n}\to\algA^{\tensor n}$ as
\begin{equation}\label{Z}
  \Z{} = \bb'\sss.
\end{equation}
They are to be viewed as the (quantum) $\overwrite{n}{Z}$ mappings in
Sec.~\ref{sec:reducible}.  Explicitly,
\begin{align*}
  \Z{2} a\tensor b&=a\tensor b - 1\tensor ab,\\
  \Z{3} a\tensor b\tensor c&=a\tensor b\tensor c
  - 1\tensor ab\tensor c + 1\tensor a\tensor bc,
\end{align*}
and so on.  For future use, we also introduce the notation
\begin{equation}\label{Z-partial}
  \Z[1,i]{j}=
  \PP_{1,i}{}\circ(\idd{(i-1)}\tensor\Z{j})
  :\algA^{\tensor(i+j-1)}\to\algA^{\tensor(i+j-1)}.
\end{equation}

Obviously, $\mult\circ\Z{2}=0$; in fact, there is the
\textit{exact} sequence
\begin{equation*}
  0\to\algA\tensor\bA\xrightarrow{\Z{2}}
  \algA\tensor\algA\xrightarrow{\mult}\algA.
\end{equation*}
To obtain a similar vanishing statement for higher $\Z{n}$, we compose
$\Z{n}$ and $\Z{n-1}$ by ``padding'' the lower one from the left with
the identity mapping:
\begin{lemma}\label{lemma:ZZ}
  For $n\geq2$,
  \begin{equation*}
    \Z{n}\circ\left(\id_1\tensor\Z{n-1}\right)=
    \left(\id_1\tensor\Z{n-1}\right)\circ{}\Z{n}=
    \bigl(\Z{2}\tensor{}
    \id^{\tensor(n-2)}\bigr)
    \circ
    \left(\id_1\tensor\Z{n-1}\right).
  \end{equation*}
  Hence, the through mappings
  \begin{gather*}
    \algA^{\tensor n}\xrightarrow{\id\tensor\Z{n-1}}
    \algA^{\tensor n}\xrightarrow{\Z{n}}\algA^{\tensor n}
    \xrightarrow{\mult_{12}}\algA^{\tensor(n-1)}\\[-2pt]
    \intertext{and}
    \algA^{\tensor n}\xrightarrow{\Z{n}}
    \algA^{\tensor n}\xrightarrow{\id\tensor\Z{n-1}}\algA^{\tensor n}
    \xrightarrow{\mult_{12}}\algA^{\tensor(n-1)}
  \end{gather*}
  are identically zero.
\end{lemma}

This lemma explains the gauge-theory interpretation of $\Z{n}$ as the
``zero modes'' $\overwrite{n}{Z}$ in Sec.~\ref{sec:reducible}.  Let
$t_a$ be a basis in~$\algA$ and
\begin{equation}\label{fabc}
  t_a\,t_b=\sum_{c}f_{ab}^c\,t_c
\end{equation}
be the algebra multiplication table.  Then the analogue of
Eqs.~\eqref{eq:involution} is obviously given by
\begin{equation}\label{fabc-skew}
  [t_a,t_b]=\sum_c\,f_{[ab]}^c\,t_c.
\end{equation}
Next, for elements $a=\sum_{c}t_{c}a^{c}$, we write
\begin{align*}
  \Z{n}(a_1\tensordots
  a_n)={}&Z_{b_1\dots b_n}^{c_1\dots c_{n-1}}
  \tensor t_{c_1}\tensordots t_{c_{n-1}}\,
  a_1^{b_1}\dots a_n^{b_n},\\
  &Z_{b_1\dots b_n}^{c_1\dots c_{n-1}}\in\algA,
\end{align*}
where summation over repeated indices is understood.  The first tensor
factor in the image of $\Z{n}$ is therefore separated from the others,
with each $Z_{b_1\dots b_n}^{c_1\dots c_{n-1}}$ being an element in
this copy of~$\algA$ (consistently with the language of noncommutative
differential forms used in what follows).

In this component form, $\Z{2}$ is represented by
$Z^a_{bc}=t_b\delta^a_c-f^a_{bc}$ and the above relation
$\mult\circ\Z{2}=0$ becomes
\begin{equation}\label{Z-zero}
  \sum_{c}Z^c_{ab}\,t_c=0,
\end{equation}
which is to be viewed as the ``reducibility'' equation~\eqref{ZT}.
Further, the ``reducibility relations'' for the higher $\Z{n}$
mappings\,---\,analogues of~\eqref{ZZ}\,---\,are the vanishing
conditions in~\bref{lemma:ZZ}, which in the component notation become
\begin{gather}\label{higher-Z-zero}
  \sum_{a_1,\dots,a_n}
  Z^{a_1\dots a_n}_{b_1\dots b_{n+1}} Z^{c_1\dots c_{n-1}}_{a_1\dots
  a_n}=0
\end{gather}
(this involves \textit{multiplication}
in~$\algA$).\footnote{\label{foot:inices}We thus see that in
  translating from the notation in Sec.~\ref{sec:reducible}, we have
  $\{\alpha_{n-1}\}=\{a_1,\dots,a_n\}$, i.e., the indices $\alpha_n$
  of each reducibility stage label elements of a basis in the $n+1$th
  tensor power of a fixed vector space.  The coefficients in the BRST
  differential are then written as, e.g.,
  $U_{\alpha_0\beta_1}^{\alpha_1}=U_{a,bc}^{de}$.  In
  Eq.~\eqref{eq:THREE}, the coefficients
  $U_{\alpha_0\alpha_1}^{\gamma_0\delta_0}=U_{a,bc}^{d,e}$ can then be
  absorbed into $U_{a,bc}^{de}$, which can be rephrased as setting
  $\Pi_{\beta_0\gamma_0}^{\alpha_0}=0$ in~\eqref{eq:THREE}
  and~\eqref{eq:SEVEN}, and similar further simplifications.}  More
precisely, Eqs.~\eqref{higher-Z-zero} are interpreted as
``\textit{quantum}'' equations~\eqref{ZZ-q}, in the particular case
where $\Pi=0$ and no $\hbar\AQ$ terms arise in the right-hand side
(consistently with the fact that $i\hbar$ is ``genuinely'' equal
to~$1$).  The component representation of the $\Z{n}$ mappings is
easily obtained from
\begin{align*}
  Z^{a_1\dots a_n}_{b_1\dots b_{n+1}} ={}& Z^{a_1\dots
    a_{n-1}}_{b_1\dots b_n}\delta^{a_n}_{b_{n+1}}
  +(-1)^n\delta^{a_1}_{b_1}\dots\delta^{a_{n-1}}_{b_{n-1}} f^{a_n}_{b_n
    b_{n+1}}.
\end{align*}

\subsection{Noncommutative differential forms}
\addcontentsline{toc}{section}{\numberline{\thesubsection.}Noncommutative
  differential forms} We recall the interpretation of
$\Omega^n=\algA\tensor \bA^{\tensor n}$ as noncommutative differential
forms~\cite{Arveson,Karoubi,[CQ]}.  The algebra of noncommutative
differential forms $\Omega^{\bullet}=\Omega^{\bullet}\algA$ over
$\algA$ is the universal differential graded algebra generated by
$\algA$ and the symbols $da$, $a\in\algA$, such that $da$ is linear in
$a$, the Leibnitz rule $d(ab) = d(a) b + adb$ is satisfied, and $d
1=0$.  The isomorphism $\algA\tensor \bA^{\tensor n}\to\Omega^n$ is
given~by
\begin{equation*}
  a_0\tensor a_1\tensordots a_n\mapsto a_0 da_1\dots da_n.
\end{equation*}
Under this isomorphism, the action of~$d$ becomes $d(a_0\tensor
a_1\tensordots a_n)= 1\tensor a_0\tensor a_1\tensordots a_n$.
Noncommutative differential forms are a bimodule over~$\algA$.  The
left action is obvious, and the rule to define the right action is as
follows: on $1$-forms, an element $c\in\algA$ acts as $a_0\,da_1\ldot
c = a_0\,d(a_1 c) - a_0a_1\,dc$, and similarly for higher-degree
forms, starting from the right end and ``propagating'' to the left
using the Leibnitz rule until all terms take the form $b_0db_1\dots
db_n$, for example
\begin{equation*}
  a_0\,da_1 da_2\ldot c = a_0\,da_1 d(a_2 c)
  - a_0\,d(a_1 a_2) d c
  + a_0 a_1\, da_2 d c.
\end{equation*}
For $D\in\Hom(\algA,\algA)$, we let $\Lie_D$ denote the Lie derivative
acting on noncommutative differential forms as~\cite{[CQ]}
\begin{multline*}
  \Lie_D(a_0\,da_1\dots da_n)
  = D(a_0)\,da_1\dots da_n\\*
  + a_0\, dD(a_1)\,da_2\dots da_n +\dots
  +a_0\,da_1\dots da_{n-1}\,dD(a_n)
\end{multline*}
and also let $\ii_D$ be the contraction
\begin{multline*}
  \ii_D(a_0\,da_1\dots da_n)
  = a_0\,D(a_1)\,da_2\dots da_n\\*
  - a_0\,da_1\ldot D(a_2)\,da_3\dots da_n
  +\dots+(-1)^{n-1}a_0\,da_1\dots da_{n-1}\ldot D(a_n).
\end{multline*}
The Cartan homotopy formula then holds,
\begin{equation*}
  \Lie_D=\ii_D d + d\ii_D.
\end{equation*}
For $a\in\algA$ and the inner derivation $[a,{-}]$ of~$\algA$, we
abuse the notation by writing $\Lie_a=\Lie_{[a,{-}]}$ and
$\ii_a=\ii_{[a,{-}]}$.

\begin{lemma} \label{lemma-commutators}
  If $D\in\Hom(\algA,\algA)$ is a derivation, then
  \begin{equation*}
    \ii_D\,\bb' + \bb'\,\ii_D{}=0
    \qquad\text{and}\qquad
    \bb'\,\Lie_D - \Lie_D\,\bb'{}=0.
  \end{equation*}
\end{lemma}

The first relation is verified directly.  Combined with the above
homotopy formula, it then immediately implies the second relation.  In
what follows, this lemma is used for inner derivations $D=[a,{-}]$; in
particular, we have $\bb'\,\Lie_a = \Lie_a\,\bb'$ for $a\in\algA$.

We extend the above relations to the bar
resolution~$\Omega^{\bullet}\tensor\algA$.

\subsection{Ghosts and the BRST differential $\BRST$}
\addcontentsline{toc}{section}{\numberline{\thesubsection.}Ghosts and
  the BRST differential~$\BRST$}
\label{sec:BRST}
We now view the bar resolution and the mappings $\Z{n}$ as a reducible
gauge theory with $t_a$ (a chosen basis in~$\algA$) playing the role
of constraints and with the reducibility relations given
by~\eqref{Z-zero} and~\eqref{higher-Z-zero}.

\subsubsection{The ghost content} \label{sec:GHOSTS}
We introduce ghosts for each term in the bar resolution,
\begin{equation*}
  \begin{split}
    &C^n\in\Omega^{n-2}\tensor\algA,\quad n\geq2,\\
    &{C}_1\in\algA
  \end{split}  
\end{equation*}
and also introduce the conjugate momenta\footnote{Relation with the
  notation for the conjugate momenta used in Sec.~\ref{sec:reducible}
  and in the literature on constrained Hamiltonian systems in general
  is~$(P_n)_{a_1,\dots,a_n}=(-1)^{n-1}\ghP_{\alpha_{n-1}}$, see
  footnote~\ref{foot:inices}.}
\begin{equation*}
  \begin{split}
    &P_n\in\Hom(\bA^{\tensor(n-1)}\tensor\algA,\oC),\quad n\geq2,\\
    &P_1\in\Hom(\algA,\oC)
  \end{split}
\end{equation*}
with the $\oZ_2$ gradings $\Grass{C^n}=\varepsilon(P_n){}\equiv{}
n\;\mathrm{mod}\,2$.  The ghost number assignments are as follows:
\begin{equation*}
  \gh C^n = n,\qquad
  \gh P_n = -n.
\end{equation*}
The $\oZ_2$ grading therefore coincides with the ghost-number grading
considered modulo~$2$; we nevertheless explicitly specify the
$\oZ_2$-grading along with the ghost-number grading.

With the elements $a_i\in\algA$ assigned $\oZ_2$-gradings
$\Grass{a_i}$, we set $\Grass{a_1\tensordots
  a_n}=\Grass{a_1}+\dots+\Grass{a_n}$ and often write $\Grass{a}$
for~$\Grass{a_1\tensordots a_n}$.

We write the canonical coupling as
\begin{equation}\label{Contr}
  \begin{split}
    (\bA^{\tensor n}\tensor\algA)
    \tensor
    \Hom(\bA^{\tensor n}\tensor\algA,\oC)
    &{}\to\oC,\\[-2pt]
    (a_1\tensordots a_n\tensor a')\tensor D&{}\mapsto
    \Contr{a_1\tensordots a_n\tensor a'}{D}.
  \end{split}
\end{equation}
This also induces the contraction
\begin{gather}\label{contr}
  \begin{split}
    (\algA\tensor\bA^{\tensor n}\tensor\algA)
    {}\tensor{}
    \Hom(\bA^{\tensor n}\tensor\algA,\oC)\to{}&\algA\\[-2pt]
    (a_0\tensordots a_n\tensor a', D)\mapsto{}&
    \contr{a_0\tensordots a_n\tensor a'}{D},
  \end{split}
\end{gather}
where
\begin{equation*}
  \contr{a_0\tensordots a_n\tensor a'}{D}
  =a_0 \Contr{a_1\tensordots a_n}{D}.
\end{equation*}

\subsubsection{Differential operators on the bar resolution and $\BRST$}
\label{subsec:Omega}
The Hamiltonian quantization prescription involves canonical
quantization of the ghosts~\cite{[BF-Ham]}; in our case,this amounts
to considering \textit{differential operators} on the bar resolution.
With each pair $C^n$, $P_n$ subject to the canonical commutation
relations (mnemonically, $[C^n,P_m]=\delta^n_m$, and hence,
$[P_n,C^m]=(-1)^{n+1}\delta^m_n$), commuting $C^n$ through $P_n$
evaluated on $a_1\tensordots a_n$ in accordance with~\eqref{Contr}
gives
\begin{multline}\label{ccr}
  C^{n}\,\langle a_1\tensordots a_n, P_n\rangle
  =(-1)^{n\Grass{a}} a_1\tensordots a_n\\*
  {}+ (-1)^{n(\Grass{a}+1)}\langle a_1\tensordots a_n, P_n\rangle C^n,
\end{multline}
or in other words,
\begin{equation*}
  [C^{n},\langle a_1\tensordots a_n, P_n\rangle]
  =(-1)^{n\Grass{a}} a_1\tensordots a_n.
\end{equation*}
We similarly evaluate commutators involving the contraction
in~\eqref{contr}, e.g.,
\begin{equation}\label{ccrii}
  [C^1, \contr{a_0\tensor a_1}{P_1}]=
  [C^1,a_0]\langle a_1, P_1\rangle
  +(-1)^{\Grass{a_0}+\Grass{a_1}} a_0 a_1,
\end{equation}
where the first term in the right-hand side involves a commutator
in~$\algA$ and the second term involves multiplication in~$\algA$.

We use the experience with the BRST formalism to seek the differential
$\BRST$ in the form
\begin{equation} \label{the-Omega}
  \BRST = \BRST_0 + \BRST_{\algA},\qquad\gh\BRST=1,
\end{equation}
where
\begin{equation}\label{Omega-A}
  \BRST_{\algA}
  = C^1 + \sum_{n\geq1}\contr{\Z{n+1} C^{n+1}}{P_n}
\end{equation}
includes the ``boundary'' terms explicitly written in~\eqref{bc}
(although it pertains to the infinitely reducible case~$\ell=\infty$)
and where
\begin{multline}\label{Omega-0}
  \BRST_0 = \mathop{\sum\!\sum}_{j\geq i\geq1}
  \Contr{\uu_{i,j}^{i+j-1}(C^i\tensor C^j)}{P_{i+j-1}}\\*[-2pt]
  {}+ \mathop{\sum\!\sum}_{j\geq i\geq1}
  \sum_{m=1}^{\intpart{\frac{i+j-2}{2}}}
  \Contr{\uu_{i,j}^{m,i+j-m-1}(C^i\tensor C^j)}{
    P_m\tensor P_{i+j-m-1}}
\end{multline}
with the ``coefficients'' $\uu$ to be determined.  This is a BRST
differential for a closed-algebra (reducible) gauge theory (see the
remarks below for its special properties).  We also write
$\BRST_{\algA}$ in the component notation
(see~\eqref{fabc}--\eqref{Z-zero}),
\begin{equation*}
  \BRST_{\algA}
  = t_a (C^1)^a + \Z{2}_{ab}^c C^{2}{}^{ab}P_1{}_{c}
  + \Z{3}_{abc}^{de} C^{3}{}^{abc}P_2{}_{de} + \dots.
\end{equation*}

\begin{Thm} \label{thm:Omega2}
\addcontentsline{toc}{section}{\numberline{\thesubsection.}Theorem}
  There exist mappings
  \begin{alignat*}{2}
    \uu_{m,m'}^{n,n'}\in{}&
    \Hom(\algA^{\tensor m}
    \tensor\algA^{\tensor m'},
    \algA^{\tensor n}\tensor\algA^{\tensor n'}),
    \quad&&
    \begin{aligned}
      &m+m'=n+n'+1,\\[-6pt]
      &1\leq n<n',\quad 1\leq m\leq m',
    \end{aligned}\\
    \uu_{m,m'}^{m+m'-1}\in{}&
    \Hom(\algA^{\tensor m}
    \tensor\algA^{\tensor m'},
    \algA^{\tensor(m+m'-1)}),
    &&1\leq m\leq m',
  \end{alignat*}
  such that the operator $\BRST$ in~\eqref{the-Omega}--\eqref{Omega-0}
  satisfies
  \begin{equation} \label{Omega2}
    \BRST\,\BRST=0.
  \end{equation}
\end{Thm}

\subsubsection*{Remarks}
\subsubsection{} The mappings $\uu_{m,m'}^{n,n'}$ whose labels do not
satisfy the restrictions in the Theorem can be considered vanishing.
It is a matter of convention that $\uu_{m,m'}^{n,n'}=0$ for $n>n'$ (an
alternative would be to impose graded symmetry with respect to the
transposition of indices), but the condition $\uu_{m,m'}^{n,n}=0$ is
essential: it implies that no ghost momentum is squared in the BRST
differential.  Another conventional condition is that
$\uu_{m,m'}^{n,n'}=0$ for~$m>m'$.  Because we only have
$\uu_{m,n}^{m',n'}$ with $m+n=m'+n'+1$, one of the four labels is
redundant, but the notation is more transparent when all labels are
kept.

\subsubsection{} We note that the operator $\BRST_{\algA}$
in~\eqref{Omega-A} is linear in both ghosts and momenta, with the
coefficients given by $\Z{n}$ mappings~\eqref{Z}; it is therefore a
differential operator with coefficients in~$\algA$, cf.~the text after
Lemma~\bref{lemma:ZZ} (and is nothing but the ``boundary terms,''
cf.~\eqref{bc}).  On the other hand, $\BRST_0$ is bilinear in the
ghosts and (separately) in the momenta but is an operator with scalar
coefficients because of the $\Contr{{-}}{P}$ and $\Contr{{-}}{P\tensor
  P}$ contractions in it.  The general structure of~$\BRST$ is the
same as described in~\eqref{gen-reducible} with additional properties
that the $\tau_A$ are linear and the $U_{AB}$ are at most bilinear in
the ghost momenta.  Conditions~\eqref{scalar-case} are satisfied
because $\tau_A$ and $U_{AB}$ are functions of only the momenta and
$U_{AB}$ are scalar as noted above.

\subsubsection{}
For the differential $\BRST$ in Theorem~\bref{thm:Omega2}, \ $\BRST_0$
is \textit{not} a differential.

\subsection{Solution for $\uu_{m,m'}^{n,n'}$ and~$\uu_{m,m'}^{m+m'-1}$
  in Theorem~\bref{thm:Omega2}} \label{sec:proof}
\addcontentsline{toc}{section}{\numberline{\thesubsection.}Solution
  for the mappings in Theorem~\ref{thm:Omega2}} 

We determine the $\uu$ mappings by solving a part of the equations
following from $\BRST\,\BRST=0$ and then show that the remaining
equations are also satisfied.  We begin with evaluating~$\BRST\,\BRST$
for the operator in~\eqref{the-Omega}--\eqref{Omega-0}.

\subsubsection{Calculating $\BRST^2$} \label{sec:proof-low} A simple
calculation using~\eqref{ccr}--\eqref{ccrii} shows that
\begin{equation*}
  \BRST_{\algA}\,\BRST_{\algA}
  = C^1C^1
  + \sum_{m\geq1}
  \contr{\Zterm_m}{P_m}
  + \mathop{\sum\!\sum}_{n\geq m\geq1}
  \contr{\Zterm_{mn}}{P_m\tensor P_n},
\end{equation*}
where the first term in the right-hand side involves the algebra
multiplication and the terms in the series are given by
$\contr{{-}}{P}$-contractions~of
\begin{align*}
  \Zterm_m={}&\ad{1}{C^1}\ldot\Z{m+1}(C^{m+1}),\\
  \Zterm_{nn}={}&
  \zeta_{n}(C^{n+1}\tensor C^{n+1}),
  \\
  \Zterm_{mn}={}&(-1)^{mn + m}
  (\mult_{1,m+2} - \mult_{m+2,1}^*)(\Z{m+1} C^{m+1}
  \tensor\Z{n+1}C^{n+1}),\quad m<n,
\end{align*}
with $\zeta_n$ given by
\begin{multline*}
  \zeta_{n}(a_1\tensordots a_{n+1}\tensor b_1\tensordots b_{n+1})=\\
  = \tfrac{1}{4}\bigl(
  [a_1,b_1]\tensor a_2\tensordots a_{n+1}\tensor b_2
  \tensordots b_{n+1}\\*
  {}-(-1)^n [a_1,b_1]\tensor b_2\tensordots b_{n+1}\tensor a_2
  \tensordots a_{n+1}\bigr).
\end{multline*}

Further using~\eqref{ccr}--\eqref{ccrii}, we next obtain
\begin{equation*}
  [\BRST_0,\BRST_{\algA}]=
  \uu_{1,1}^1(C^1\tensor C^1)
  + \sum_{m\geq1}
  \contr{\Yterm_m}{P_m}
  + \mathop{\sum\!\sum}_{n\geq m\geq1}
  \contr{\Yterm_{mn}}{P_m\tensor P_n},
\end{equation*}
where the terms in the series are $\contr{{-}}{P}$-contractions of
\begin{multline*}
  \Yterm_{mn}=  
  (-1)^{n} \Z[1,m+1]{n+1}
  \circ
  \mathop{\sum_{j\geq i\geq0}}_{i+j=m+n}
  \uu_{i+1,j+1}^{m,n+1}(C^{i+1}\tensor C^{j+1})\\
  {}+ \mathop{\sum_{j>i\geq0}}_{i+j=m+n}
  (-1)^{m+n}
  (\id\tensor\uu_{i+1,j}^{m,n})
  \circ\Z[1,i+2]{j+1}(C^{i+1}\tensor C^{j+1})+{}
  \\
  {}+\mathop{\sum_{j+1\geq i\geq1}}_{i+j=m+n}
  (-1)^{i+1}
  (\id\tensor\uu_{i,j+1}^{m,n})
  (\Z{i+1}C^{i+1}\tensor C^{j+1})\\
  {}+ (-1)^{m+n}
  (\Z{m+1}\tensor\idd{n})
  \,\circ\!\!
  \mathop{\sum_{j\geq i\geq0}}_{i+j=m+n}\!\!
  \uu_{i+1,j+1}^{m+1,n}(C^{i+1}\tensor C^{j+1}),~m<n,
\end{multline*}
and
\begin{equation*}
  \Yterm_{nn}=  
  (-1)^{n}
  \Sym{n}{2,n+1;n+2,2n+1}\bigl(
  \Z[1,n+1]{n+1}
  \circ
  \sum_{i=0}^{n}
  \uu_{i+1,2n-i+1}^{n,n+1}(C^{i+1}\tensor C^{2n-i+1})\bigr).
\end{equation*}
We here use the notation~\eqref{Z-partial}, and the operator
$\Sym{n}{p,q;r,s}$, with $n$ considered modulo~$2$ and $p\leq q<r\leq
s$ such that $s-r=q-p$, performs graded symmetrization (for $n$ even)
or antisymmetrization (for $n$ odd) of tensor factors in the positions
$[p,\dots,q]$ and $[r,\dots,s]$, for example,
\begin{multline*}
  \Sym{n}{2,3;4,5}(a\tensor b\tensor c\tensor d\tensor e) =
  \thalf\,a\tensor b\tensor c\tensor d\tensor e\\*
  {}+ (-1)^n\cdot(-1)^{(\Grass{b}+\Grass{c})(\Grass{d}+\Grass{e})}\,
  \thalf\,a\tensor d\tensor e\tensor b\tensor c.
\end{multline*}
We also find that $\Yterm_m=\Yterm_{0,m}$ if we set
\begin{equation} \label{U0}
  \uu_{i,j}^{0,\ell}=\uu_{i,j}^{\ell}.
\end{equation}

Finally, calculating $\BRST_0\BRST_0$ gives
\begin{gather*}
  \BRST_0\BRST_0=
  \sum_{m\geq2}
  \Contr{\Xterm_m}{P_m}
  + \smash[b]{\mathop{\sum_{m\geq1}\!\sum_{n\geq1}}_{m\leq n}}
  \Contr{\Xterm_{mn}}{P_m\tensor P_n}
  + \Xterm^{(3)},
\end{gather*}
where
\begin{gather*}
  \Xterm_{mn}=
  {
    \mathop{\sum_{i\geq1}\sum_{j\geq1}}_{\substack{2i\leq m+n+1\\[1pt]
        2j\leq m+n+2}}}
  (-1)^{(i+1)(m+n)+1}
  \uu_{i,m+n+1-i}^{m,n}
  \circ
  \uu_{j,m+n+2-j}^{i,m+n+1-i}(C^{j}\tensor C^{m+n+2-j}),
\end{gather*}
$\Xterm_m=\Xterm_{0m}$, and $\Xterm^{(3)}$ denotes third-order terms
in the $C^n$ ghosts.  The square of $\BRST$ is therefore given by
\begin{multline}\label{summarize}
  \BRST\,\BRST
  = \Zterm_0
  + \sum_{m\geq1}
  \contr{\Zterm_m + \Yterm_m + 1\tensor\Xterm_m}{P_m}\\*
  + \mathop{\sum\!\sum}_{n\geq m\geq1}
  \contr{\Zterm_{mn} + \Yterm_{mn} + 1\tensor\Xterm_{mn}}{
    P_m\tensor P_n}
  + \Xterm^{(3)},
\end{multline}
where $\Zterm_0=C^1C^1 + \uu_{1,1}^1(C^1\tensor C^1)$ and
$\Zterm_{m}$, $\Yterm_{m}$, $\Xterm_{m}$, $\Zterm_{mn}$,
$\Yterm_{mn}$, $\Xterm_{mn}$ with $n\geq m\geq1$ are given above (as
we have seen, $\Xterm_1=0$ and also $\Xterm_{nn}=0$).
Expanding~\eqref{summarize} in the ghosts and momenta and equating
each power to zero, we obtain a specific form of
Eqs.~\eqref{C1}--\eqref{C3}.  In particular, the equations that are
cubic in~$C^n$, $\Xterm^{(3)}=0$, are nothing but a rewriting
of~\eqref{C3}.  The other terms in~\eqref{summarize} are at most
quadratic in the $C^n$ ghosts.  Their dependence on the ghost momenta
is shown in~\eqref{summarize} explicitly, and we therefore first
expand in the momenta.  The equation $\BRST\,\BRST=0$ is then
equivalent to the set of equations $\Xterm^{(3)}=0$ and
\begin{equation} \label{the-equations}
  \begin{split}
    \Zterm_m + \Yterm_m + 1\tensor\Xterm_m &= 0,
    \quad 1\leq m,\\
    \Zterm_{mn} + \Yterm_{mn} + 1\tensor\Xterm_{mn} &= 0,    
    \quad 1\leq m<n,\\
    \Sym{n}{2,n+1;n+2,2n+1}
    (\Zterm_{nn} + \Yterm_{nn} + 1\tensor\Xterm_{nn}) &= 0,
    \quad 1\leq n.
  \end{split}
\end{equation}
Further expanding each of these in the $C^n$ ghosts gives an infinite
list of equations starting with~\eqref{eq:THREE}--\eqref{eq:SIX}.
Speaking about the equations $\Zterm_{mn} + \Yterm_{mn} +
1\tensor\Xterm_{mn}\!=\!0$ in general, we often mean the
$\Sym{n}{2,n+1;n+2,2n+1}$ symmetrization of the ($m\!=\!n$) equations
without explicitly specifying it in the notation.  We now solve
Eqs.~\eqref{the-equations} by taking a certain projection of these
equations that yields a set of recursive relations for the sought
mappings~$\uu_{*,*}^{**}$.

\subsubsection{Finding the lowest mappings}  The lowest-order terms
$\Zterm_0$ in~\eqref{summarize} vanish if $\uu_{1,1}^1=-\mult$, minus
the multiplication in~$\algA$.  Because $C^1$ is $\oZ_2$-odd, only the
antisymmetric part of the multiplication actually contributes, and we
have
\begin{equation}\label{u111}
  \uu_{1,1}^1=-\thalf[{-},{-}]
\end{equation}
(which is a totally standard BRST fact that the lowest coefficient in
the BRST operator involves Lie algebra structure constants).  In the
order $P_1$, we readily find that
\begin{equation*}
  \Zterm_1+\Yterm_1
  = \Lie_{C^1}{}\ldot{}\Z{2}C^2 - \Z{2}\circ\uu_{1,2}^{2}
  (C^1\tensor C^2)
\end{equation*}
(the term $1\tensor\Xterm_1$ vanishes).  But
\begin{equation*}
  \Lie_a{}\circ{}\Z{2} = \Lie_a \bb'\,d
  \stackrel{\bref{lemma-commutators}}{=}\bb'\Lie_a d
  = \bb'(d\ii_a)d=\bb' d(\ii_a d + d \ii_a)=\Z{2}\circ\Lie_a,
\end{equation*}
which shows that $\Zterm_1+\Yterm_1=0$ if
\begin{equation} \label{u122}
  \uu_{1,2}^2(C^1\tensor C^2)=\Lie_{C^1}C^2.
\end{equation}
We note that this solves Eq.~\eqref{eq:THREE} in the list of equations
in~\bref{sec:eq-list}.

\subsubsection{Establishing a recursion} \label{sec:recursion} With
the lowest two mappings $\uu_{1,1}^1$ and $\uu_{1,2}^2$ thus
determined, we write the $\Z{n}$ mappings as $\Z{}=\id - 1\tensor
\bb'$ and consider the terms that do \textit{not} have the form
$1\tensor\ldots\,$ in each equation in~\eqref{the-equations}.  The
operation $\mathcal{P}_{\not{\kern1pt1}}$ of projecting onto such
terms amounts to dropping all explicit occurrences of $1\tensor\dots$
(in particular, of $1\tensor\Xterm_{mn}$) and replacing
$\Z{n}\to\idd{n}$ and $\Z[1,i]{n}\to\PP_{1,i}\rule[-12pt]{0pt}{12pt}$.
Equations~\eqref{the-equations} thus imply the equations
\begin{multline*}\tag{$\boldsymbol{m,\!n}$}
  \sum_{i=0}^{\intpart{\frac{m+n}{2}}}
  \uu_{i+1,m+n-i+1}^{m,n+1}(C^{i+1}\tensor C^{m+n-i+1})=\\
  \shoveleft{{}=-(-1)^{mn + m + n}
  \PP_{m+1,1}{}\circ{}
  (\mult_{1,m+2} - \mult_{m+2,1}^*)(C^{m+1}\tensor C^{n+1})}\\*
  \qquad{}- (-1)^{m}\sum_{i=0}^{\intpart{\frac{m+n-1}{2}}}
  \PP_{m+1,1}{}\circ{}(\id\tensor\uu_{i+1,m+n-i}^{m,n})
  \circ\PP_{1,i+2}(C^{i+1}\tensor C^{m+n-i+1})\\*
  \qquad{}+\sum_{i=1}^{\intpart{\frac{m+n+1}{2}}}
  (-1)^{i+n}
  \PP_{m+1,1}{}\circ{}(\id\tensor\uu_{i,m+n-i+1}^{m,n})
  (C^{i+1}\tensor C^{m+n-i+1})\\*[-4pt]
  {}-(-1)^{m}
  \sum_{i=0}^{\intpart{\frac{m+n}{2}}}
  \PP_{m+1,1}\circ\uu_{i+1,m+n-i+1}^{m+1,n}(C^{i+1}\tensor C^{m+n-i+1}),
\end{multline*}
where $0\leq m\leq n$ (and we recall Eq.~\eqref{U0} for $m=0$).  Each
of these relations amounts to $\intpart{\frac{m+n}{2}}+1-m$
independent equations obtained by extracting each ghost pair
$(C^a,C^b)$, but to save space, we keep them in the above form of
``generating relations''.  We also remark that in applying mappings
to $C^\ell\tensor C^\ell$, graded symmetrization with respect to the
two ``halves'' of the tensor argument must also be performed in
accordance with~$\Grass{C^\ell}=\ell$.

To see that Eqs.~\thetag{$\boldsymbol{m,\!n}$} are in fact a set of
recursive relations, we fix $m+n=N$ with a positive integer~$N$ and
arrange the $\half(\intpart{\frac{N}{2}}+1)(\intpart{\frac{N}{2}}+2)$
equations following from~\thetag{$\boldsymbol{m,\!n}$} in the order
specified by
\begin{equation}\label{the-order}
  \begin{aligned}
  (m,n)={}&(\intpart{\tfrac{N}{2}}, \intpart{\tfrac{N+1}{2}})
  &&(1\text{ equation}),\\
  (m,n)={}&(\intpart{\tfrac{N}{2}}-1, \intpart{\tfrac{N+1}{2}}+1)
  \qquad&&(2\text{ equations}),\\
  \ldots\\
  (m,n)={}&(1, N-1)
  &&(\intpart{\tfrac{N}{2}}\text{ equations}),\\
  (m,n)={}&(0,N)
  &&(\intpart{\tfrac{N}{2}}+1\text{ equations}).
\end{aligned}
\end{equation}
The underlying partial ordering is therefore given by
\begin{equation}\label{partial-ordering}
  \uu_{m,n}^{m',n'}\prec\uu_{m_1,n_1}^{m'_1,n'_1}
  \text{ if }
  \begin{array}[t]{l}
    \text{either } m'+n'<m'_1+n'_1\\
    \text{or } m'+n'=m'_1+n'_1 \text{ and } n'-m'<n'_1-m'_1.
  \end{array}
\end{equation}
It follows that for each integer $N\geq4$, the vanishing
$\uu_{m,n}^{m',n'}$ mappings are given by
\begin{equation}
  \uu_{\intpart{\frac{N}{2}}-i,\intpart{\frac{N+1}{2}}+i}^{
    \intpart{\frac{N}{2}}-j,\intpart{\frac{N+1}{2}}+j-1} = 0,\quad
  1\leq j\leq i\leq\intpart{\tfrac{N}{2}}-1.
\end{equation}
In other words, 
\begin{equation} \label{eq:vanishing}
  \uu_{m,n}^{m',m+n-m'-1}=0\quad\text{if}\quad m'\geq m
  \quad\text{and}\quad n-m\geq2.
\end{equation}
The recursive relations are now written most conveniently for odd and
even $m+n$ separately.

\subsubsection{Odd $\boldsymbol{m+n}$} For $m+n=2k-1$,
Eqs.~\thetag{$\boldsymbol{m,\!n}$}, which can then be labeled
as~\thetag{$\boldsymbol{2k-1;m}$} with $0\leq m\leq k-1$, become
\begin{multline*}\tag{$\boldsymbol{2k\!-\!1;m}$}
  \sum_{i=m}^{k-1}
  \uu_{i+1,2k-i}^{m,2k-m}(C^{i+1}\tensor C^{2k-i})=\\
  \shoveleft{\quad{}=\PP_{m+1,1}{}\circ{}
    \Bigl((\mult_{1,m+2} - \mult_{m+2,1}^*)(C^{m+1}\tensor C^{2k-m})}\\*
  {}- (-1)^{m}\sum_{i=m}^{k-1}
  (\id\tensor\uu_{i+1,2k-i-1}^{m,2k-m-1})
  \circ\PP_{1,i+2}(C^{i+1}\tensor C^{2k-i})\\*
  {}-\sum_{i=m+1}^{k}(-1)^{i+m}
  (\id\tensor\uu_{i,2k-i}^{m,2k-m-1})
  (C^{i+1}\tensor C^{2k-i})\\*[-4pt]
    {}-(-1)^{m}
    \sum_{i=m+1}^{k-1}
    \uu_{i+1,2k-i}^{m+1,2k-m-1}(C^{i+1}\tensor C^{2k-i})\Bigr).
\end{multline*}
It is now obvious that these relations express
$\uu_{i+1,2k-i}^{m,2k-m}$ through
$\uu_{a,b}^{c,d}\prec\uu_{i+1,2k-i}^{m,2k-m}$.  In more detail, we
first choose the top value $m=k-1$ in accordance with the ordering.
The last sum in the right-hand side is then absent, and
$\uu_{k,k+1}^{k-1,k+1}$ is therefore expressed through
$\uu_{m,n}^{m',n'}$ with $m'+n'=2k-1<2k$, namely,
\begin{multline}\label{rec-next}
  \uu_{k,k+1}^{k-1,k+1}(C^{k}\tensor C^{k+1})
  {}=\PP_{k,1}{}\circ{}
  \Bigl((\mult_{1,k+1} - \mult_{k+1,1}^*)(C^{k}\tensor C^{k+1})\\*
  {}+ (-1)^{k}
  (\id\tensor\uu_{k,k}^{k-1,k})
  \circ\PP_{1,k+1}(C^{k}\tensor C^{k+1})
  + (\id\tensor\uu_{k,k}^{k-1,k})
  (C^{k+1}\tensor C^{k})\Bigr),
\end{multline}
Next, setting $m=k-2$ gives the two equations
\begin{equation}\label{recursion-one}
  \uu_{k-1,k+2}^{k-2,k+2}
  =\PP_{k-1,1}\circ\Bigl(\mult_{2,1}^* - \mult_{1,2}
  -(-1)^k(\id\tensor\uu_{k-1,k+1}^{k-2,k+1})
  \Bigr)\circ\PP_{1,k}
\end{equation}
and
\begin{multline}\label{rec-second}
  \uu_{k,k+1}^{k-2,k+2}(C^k\tensor C^{k+1})=
  \PP_{k-1,1}\circ\Bigl(
  -(-1)^k \uu_{k,k+1}^{k-1,k+1}(C^{k}\tensor C^{k+1})
  \\
  + (\id\tensor\uu_{k-1,k+1}^{k-2,k+1})(C^{k}\tensor C^{k+1})
  - (\id\tensor\uu_{k,k}^{k-2,k+1})(C^{k+1}\tensor C^{k})\\
  -(-1)^k (\id\tensor\uu_{k,k}^{k-2,k+1})
  \circ\PP_{1,k+1}(C^{k}\tensor C^{k+1})\Bigr),
\end{multline}
and~so~on: for each $m$, $0\leq m\leq k-2$, the equations for
``generic'' values of $i$, i.e., $m+1\leq i\leq k-1$, are given by
\begin{multline}\label{gen-recursion}
  \uu_{i+1,2k-i}^{m,2k-m}
  =\PP_{m+1,1}{}\circ{}
  \Bigl(
  (-1)^{m+1}(\id\tensor\uu_{i+1,2k-i-1}^{m,2k-m-1})
  \circ\PP_{1,i+2}\\*
  {}- \id\tensor\uu_{i,2k-i}^{m,2k-m-1}
  -(-1)^{m}\uu_{i+1,2k-i}^{m+1,2k-m-1}\Bigr).
\end{multline}
The ``boundary'' equations (those with $i=m$ and $i=k$) are also
easily extracted from~\thetag{$\boldsymbol{2k\!-\!1;m}$} (the
($i=m$)-equation involves the ``inhomogeneous'' contribution
$\PP_{m+1,1}{}\circ{}(\mult_{1,m+2} - \mult_{m+2,1}^*)$ in the
right-hand side).

\subsubsection{Even $\boldsymbol{m+n}$} For $m+n=2k$,
Eqs.~\thetag{$\boldsymbol{m,\!n}$}, now labeled
as~\thetag{$\boldsymbol{2k;m}$} with $0\leq m\leq k$, become
\begin{multline*} \tag{$\boldsymbol{2k;m}$}
  \sum_{i=m}^{k}
  \uu_{i+1,2k-i+1}^{m,2k-m+1}(C^{i+1}\tensor C^{2k-i+1})=\\
  \shoveleft{\quad{}=-(-1)^{m}
    \PP_{m+1,1}{}\circ{}
    \Bigl(
    (\mult_{1,m+2} - \mult_{m+2,1}^*)(C^{m+1}\tensor C^{2k-m+1})}\\*
  {}-(-1)^{m}\sum_{i=m}^{k-1}
  (\id\tensor\uu_{i+1,2k-i}^{m,2k-m})
  \circ\PP_{1,i+2}(C^{i+1}\tensor C^{2k-i+1})\\*
  \quad{}+\sum_{i=m+1}^{k}
  (-1)^{i+m}
  (\id\tensor\uu_{i,2k-i+1}^{m,2k-m})
  (C^{i+1}\tensor C^{2k-i+1})\\
    {}-(-1)^{m}
    \sum_{i=m+1}^{k}
    \uu_{i+1,2k-i+1}^{m+1,2k-m}(C^{i+1}\tensor C^{2k-i+1})\Bigr).
\end{multline*}
These equations must also be considered in the order specified by
consecutively taking $m=k$, $m=k-1$, \dots, $m=1$.  The top value
$m=k$ is somewhat special here: on one hand, all the
$\uu_{m,n}^{m',n'}$ mappings drop from the right-hand side and only
the first term survives, thereby giving an explicit expression for
$\uu_{k+1,k+1}^{k,k+1}$; on the other hand, the equation must only be
satisfied after the $\Sym{k}{2,k+1;k+2,2k+1}$ symmetrization (see the
remark at the end of~\bref{sec:proof-low}) and with the proper
symmetrization of the ghost argument $C^{k+1}\tensor C^{k+1}$ (see the
remark after equation~\thetag{$\boldsymbol{m,\!n}$}).  This leaves
some freedom in determining $\uu_{k+1,k+1}^{k,k+1}$, which can be
fixed using the \textit{full} equation $\Sym{k}{2,k+1;k+2,2k+1}
(\Zterm_{kk}+\Yterm_{kk}+1\tensor\Xterm_{kk})=0$ for~$\uu_{k, k}^{k-1,
  k}$.  It follows that
\begin{equation}\label{U-special2}
  \uu_{k, k}^{k-1, k}
  = (-1)^{k+1}\bb'\tensor\idd{k}{}
  + (-1)^k(\idd{(k-1)}{}\tensor{}\bb')
  \circ{\PP_{k,1}},
\end{equation}
where $\PP_{k,1}$ denotes the operator inverse to $\PP_{1,k}$ and
graded symmetrization with respect to the two ``halves'' of the tensor
argument must be performed in evaluating $\uu_{k, k}^{k-1, k}$ on
$C^k\tensor C^k$, for example (omitting sign factors due to the
$\oZ_2$-grading),
\begin{multline}\label{u2212}
  \uu_{2,2}^{1,2}\,a\tensor b\tensor a' \tensor b'
  = \thalf(-a b\tensor a'\tensor b'
  + b\tensor a a'\tensor b' - b\tensor a\tensor a' b'\\*
  -a' b'\tensor a\tensor b
  + b'\tensor a' a\tensor b - b'\tensor a'\tensor a b)
\end{multline}
(which actually solves Eq.~\eqref{eq:SEVEN}).  

The next (two) equations following from~\thetag{$\boldsymbol{2k;m}$},
which correspond to $m=k-1$, allow expressing $\uu_{k,k+2}^{k-1,k+2}$
through $\uu_{k,k+1}^{k-1,k+1}$ and $\uu_{k+1,k+1}^{k-1,k+2}$ through
$\uu_{k,k+1}^{k-1,k+1}$ and~$\uu_{k+1,k+1}^{k,k+1}$; all these
equations are easily written out similarly to~\eqref{gen-recursion}.
In applying the recursion further, we must only take into account that
for $0\leq m\leq k-1$, the equations
\begin{multline}\label{uk1k1metc}
  \uu_{k+1,k+1}^{m,2k-m+1}(C^{k+1}\tensor C^{k+1})=\\  
  {}=(-1)^{k+m}
  \PP_{m+1,1}{}\circ{}(\id\tensor\uu_{k,k+1}^{m,2k-m})
  (C^{k+1}\tensor C^{k+1})\\
  {}-(-1)^{m}
  \PP_{m+1,1}\circ
  \uu_{k+1,k+1}^{m+1,2k-m}(C^{k+1}\tensor C^{k+1}),
\end{multline}
implied by~\thetag{$\boldsymbol{2k;m}$} are evaluated on
$C^{k+1}\tensor C^{k+1}$ and the mappings in the right-hand side must
therefore be symmetrized appropriately with respect to the two
``halves'' of the tensor argument (i.e., antisymmetrized for $k+1$ odd
and symmetrized for $k+1$ even).

Equations~\eqref{u111}, \eqref{u122}, \eqref{eq:vanishing},
\thetag{$\boldsymbol{2k\!-\!1;m}$}, and~\thetag{$\boldsymbol{2k;m}$}
determine all the mappings involved in the differential~$\BRST$.

\begin{rem}
  It follows that the structure in Eq.~\eqref{u122} propagates through
  the recursive relations to
  \begin{equation}\label{u1nn}
    \uu_{1, n}^{n}\,a_0\tensor a_1\tensordots a_n
    = (-1)^n\Lie_{a_0}\, a_1\tensordots a_n.
  \end{equation}
\end{rem}

\subsubsection{Examples} In addition to~\eqref{u111}, \eqref{u122},
\eqref{u2212} (and~\eqref{u1nn}), we give explicit formulas for
several lowest mappings that can easily be derived from the recursive
relations above.  For $\uu_{m,n}^{m',n'}$ with $m'+n'=3$, we have,
along with~\eqref{u2212},
\begin{multline*}
  2\,\uu_{2, 2}^{3}\,a\tensor b\tensor a'\tensor b' =
  a\tensor a'\tensor b' b
  - a\tensor a'\tensor b b'
  + a\tensor a' b\tensor b' - a\tensor b a'\tensor b'\\*
  + a'\tensor a\tensor b b' - a'\tensor a\tensor b' b
  + a'\tensor a b'\tensor b - a'\tensor b' a\tensor b
  + b\tensor a\tensor a' b'\\*
  - b\tensor a a'\tensor b'
  + b'\tensor a'\tensor a b - b'\tensor a' a\tensor b + 
  a b\tensor a'\tensor b' + a' b'\tensor a\tensor b.
\end{multline*}
One of the $\uu_{m,n}^{m',n'}$ mappings with $m'+n'=4$ is given
by
\begin{multline*}
  \uu_{2, 3}^{1, 3}\,a\tensor b\tensor a'\tensor b'\tensor c' =
  c'\tensor a'\tensor b' a\tensor b
  - c'\tensor a'\tensor b'\tensor a b\\
  - b\tensor a' a\tensor b'\tensor c'
  + b\tensor a'\tensor a b'\tensor c'
  -b\tensor a'\tensor a\tensor b' c'\\
  + b\tensor a a'\tensor b'\tensor c'
  - a b\tensor a'\tensor b'\tensor c'
  - b' c'\tensor a'\tensor a\tensor b,
\end{multline*}
or in the component notation (see~\eqref{fabc}--\eqref{Z-zero}), by
the tensor
\begin{multline*}
  \uu_{ab,cde}^{g,hij} =
  \delta_e^g\delta_c^h f_{d a}^i\delta_b^j
  - \delta_e^g\delta_c^h\delta_d^i f_{a b}^j
  - \delta_b^g f_{c a}^h\delta_d^i\delta_e^j\\*
  {}+ \delta_b^g\delta_c^h f_{a d}^i\delta_e^j
  -\delta_b^g\delta_c^h\delta_a^i f_{d e}^j
  + \delta_b^g f_{a c}^h\delta_d^i\delta_e^j
  - f_{a b}^g\delta_c^h\delta_d^i\delta_e^j
  - f_{d e}^g\delta_c^h\delta_a^i\delta_b^j.
\end{multline*}

\subsection{The end of the proof of
  Eqs.~(\bref{the-equations})}\label{sec:end-proof-eq}
\addcontentsline{toc}{section}{\numberline{\thesubsection.}The end of
  the proof of Eqs.~\eqref{the-equations}}

With the $\uu_{m,n}^{\ell,m+n-\ell-1}$ mappings found recursively from
the $\mathcal{P}_{\not{\kern1pt1}}$-projections of
Eqs.~\eqref{the-equations}, we must next show that the remaining parts
of the equations
\begin{equation*}
  \Zterm_{mn} + \Yterm_{mn} + 1\tensor\Xterm_{mn} = 0,
  \quad 0\leq m<n
\end{equation*}
are satisfied.  Explicitly, the equations to be verified are obtained
by replacing $\Z{}\to-1\tensor\bb'$; all tensor terms then acquire the
form $1\tensor\mathbb{X}$, and the resulting equations $\mathbb{X}=0$
(as previously, written as ``generating equations'') are given by
\begin{multline}\label{rest-eqs}
  (-1)^{m} (\idd{m}\tensor\bb')\circ
  \sum_{i=m}^{\intpart{\frac{m+n}{2}}}
  \uu_{i+1,m+n-i+1}^{m,n+1}(C^{i+1}\tensor C^{m+n-i+1})\\*
  {}+ (\bb'\tensor\idd{n})
  \circ
  \sum_{i=m+1}^{\intpart{\frac{m+n}{2}}}
  \uu_{i+1,m+n-i+1}^{m+1,n}(C^{i+1}\tensor C^{m+n-i+1})\\
  {}+\sum_{i=m}^{\intpart{\frac{m+n-1}{2}}}
  \uu_{i+1,m+n-i}^{m,n}(C^{i+1}\tensor\bb' C^{m+n-i+1})\\
  {}+\sum_{i=m+1}^{\intpart{\frac{m+n+1}{2}}}
  (-1)^{m+n+i+1}
  \uu_{i,m+n-i+1}^{m,n}
  (\bb'C^{i+1}\tensor C^{m+n-i+1})\\
  {}+\!\sum_{j=m+1}^{\intpart{\frac{m+n}{2}}}
  \sum_{i=j}^{\intpart{\frac{m+n}{2}}}
  (-1)^{j(m+n)}
  \uu_{j,m+n+1-j}^{m,n}
  \circ
  \uu_{i+1,m+n+1-i}^{j,m+n+1-j}(C^{i+1}\tensor C^{m+n+1-i})=0
\end{multline}
for $0\leq m<n$.  We must show that they are satisfied with the
mappings~$\uu_{m,n}^{\ell,m+n-\ell-1}$ determined by the above
recursion.

\addtocounter{subsubsection}{-1}
\subsubsection{The strategy} Equations~\eqref{rest-eqs} are proved by
(somewhat tedious) induction on the order of the $\uu$ mappings
introduced in~\eqref{partial-ordering}.  We consider
Eqs.~\eqref{rest-eqs} with $m+n=2k-1$ and assume that all equations
with $m+n\leq2k-2$ are satisfied (Eqs.~\eqref{rest-eqs} with $m+n=2k$
are proved similarly).  The equations that we must verify for
$m+n=2k-1$ are as follows.  Setting $m=k-1-\ell$ and $n=k+\ell$, we
have $\ell+1$ equations for each $\ell\in\{0,\dots,k-1\}$; for
$\ell\geq1$, these are
\begin{multline}\label{i=0}
  (-1)^{k+\ell}(\idd{(k-\ell-1)}\tensor\bb')
  \circ\uu_{k,k+1}^{k-\ell-1,k+\ell+1}(C^k\tensor C^{k+1})\\*
  -(\bb'\tensor\idd{(k+\ell)})\circ\uu_{k,k+1}^{k-\ell,k+\ell}
  (C^k\tensor C^{k+1})
  -\uu_{k,k}^{k-\ell-1,k+\ell}(C^k\tensor\bb'C^{k+1})\\*
  +(-1)^{k}\uu_{k-1,k+1}^{k-\ell-1,k+\ell}
  (\bb'C^k\tensor C^{k+1})
  -(-1)^{k}\uu_{k,k}^{k-\ell-1,k+\ell}
  (\bb'C^{k+1}\tensor C^{k})\\*
  -\sum_{j=0}^{\ell-1}(-1)^{k+1+j}\uu_{k-j-1,k+j+1}^{k-\ell-1,k+\ell}
  \circ\uu_{k,k+1}^{k-j-1,k+j+1}(C^k\tensor C^{k+1})=0,
\end{multline}
\begin{gather}\kern-10pt
  (-1)^{k+\ell}(\idd{(k-\ell-1)}\!\tensor\!\bb')
  \circ\uu_{k-\ell,k+\ell+1}^{k-\ell-1,k+\ell+1}
  -\uu_{k-\ell,k+\ell}^{k-\ell-1,k+\ell}\circ(\idd{(k-\ell)}
  \tensor\bb')
  =0,\kern-10pt
\end{gather}
and
\begin{multline}\label{ell-ge-2}
  (-1)^{k+\ell}(\idd{(k-\ell-1)}\tensor\bb')
  \circ\uu_{k-i,k+i+1}^{k-\ell-1,k+\ell+1}
  -(\bb'\tensor\idd{(k+\ell)})\circ\uu_{k-i,k+i+1}^{k-\ell,k+\ell}\\*
  -\uu_{k-i,k+i}^{k-\ell-1,k+\ell}\circ(\idd{(k-i)}\tensor\bb')
  +(-1)^{k+i}\uu_{k-i-1,k+i+1}^{k-\ell-1,k+\ell}\circ
  (\bb'\tensor\idd{(k+1+i)})\\*
  -\sum_{j=i}^{\ell-1}(-1)^{k+1+j}\uu_{k-j-1,k+j+1}^{k-\ell-1,k+\ell}
  \circ\uu_{k-i,k+i+1}^{k-j-1,k+j+1}=0,\quad 1\leq i\leq\ell-1.
\end{multline}
For~$\ell=0$, the only equation is
\begin{multline}\label{ell=0}
  (\idd{(k-1)}\tensor\bb')\circ
  \uu_{k,k+1}^{k-1,k+1}(C^{k}\tensor C^{k+1})={}\\
  =(-1)^{k}\uu_{k,k}^{k-1,k}(C^{k}\tensor\bb' C^{k+1})
  + \uu_{k,k}^{k-1,k}
  (\bb'C^{k+1}\tensor C^{k}).
\end{multline}

For a fixed $k$, we proceed by induction on~$\ell$, which is actually
part of the induction on the order of the $\uu$ mappings
in~\eqref{partial-ordering} --- $\ell$~labels rows
in~\eqref{the-order}. At each step, we use the defining recursive
relations and the lower (previously proved)
equations~\eqref{i=0}--\eqref{ell=0}.

\subsubsection{} 
We begin with the $\ell$-induction base~$\ell=0$ and show
that~\eqref{ell=0} is satisfied.  With~\eqref{U-special2},
Eq.~\eqref{ell=0} becomes
\begin{multline}\label{work-90}
  (\idd{(k-1)}\tensor\bb')\circ
  \uu_{k,k+1}^{k-1,k+1}(C^{k}\tensor C^{k+1})={}\\
  -\bb'C^k\tensor\bb'C^{k+1}    
  {}+(\idd{(k-1)}\tensor\bb')\circ\PP_{k,1}
  \bigl(C^k\tensor\bb'C^{k+1} 
  + (-1)^k\bb'C^{k+1}\tensor C^k\bigr).
\end{multline}
Using the obvious identity
\begin{equation}\label{obvious}
  (\idd{(k-1)}\tensor\bb')\circ\PP_{k,1}=
  \PP_{k,1}\circ\bigl(\mult_{1,k+1} - 
  (\idd{k}\tensor\bb')\bigr),
\end{equation}
we further rewrite~\eqref{work-90} as
\begin{multline}\label{work-91}
  (\idd{(k-1)}\tensor\bb')\circ
  \uu_{k,k+1}^{k-1,k+1}(C^{k}\tensor C^{k+1})={}\\
  \shoveleft{\quad{}=  -\bb'C^k\tensor\bb'C^{k+1}
    -(-1)^k\PP_{k,1}(\bb'C^{k+1}\tensor\bb'C^k)}\\*
  +\PP_{k,1}\circ\mult_{1,k+1}(C^k\tensor\bb'C^{k+1})
  +(-1)^k\PP_{k,1}\circ\mult_{1,k+1}(\bb'C^{k+1}\tensor C^k).
\end{multline}
On the other hand, 
the mapping $\uu_{k,k+1}^{k-1,k+1}$ is expressed via recursive
relation~\eqref{rec-next}, and therefore (applying~\eqref{obvious}
again), we can rewrite the left-hand side of~\eqref{work-91} as
\begin{multline*}
  (\idd{(k-1)}\tensor\bb')\circ
  \uu_{k,k+1}^{k-1,k+1}(C^{k}\tensor C^{k+1})={}\\
  \shoveleft{{}=\PP_{k,1}{}\circ\mult_{1,k+1}\circ{}
    \Bigl((\mult_{1,k+1} - \mult_{k+1,1}^*)(C^{k}\tensor C^{k+1})}\\*
  {}+ (-1)^{k}
  (\id\tensor\uu_{k,k}^{k-1,k})
  \circ\PP_{1,k+1}(C^{k}\tensor C^{k+1})
  + (\id\tensor\uu_{k,k}^{k-1,k})
  (C^{k+1}\tensor C^{k})\Bigr)\\*
  -\PP_{k,1}\circ(\idd{k}{}\tensor\bb')
  \circ(\mult_{1,k+1} - \mult_{k+1,1}^*)(C^{k}\tensor C^{k+1})
  +\mult_{k.k+1}(\bb'C^k\tensor C^{k+1})\\*
  -\bb'C^k\tensor\bb'C^{k+1}
  +(-1)^k\PP_{k,1}((\id\tensor\bb')C^{k+1})\tensor\bb'C^k,
\end{multline*}
where we also used that $(\id\tensor\bb'_{(k)})\circ\PP_{1,k+1}=
\PP_{1,k}\circ(\bb'_{(k)}\tensor\id)$ (with the subscript in
$\bb'_{(k)}$ indicating that $\bb'$ acts on $k$ tensor factors).  We
now have
\begin{multline*}
  \text{left-hand side of~\eqref{work-91}} -
  \text{right-hand side of~\eqref{work-91}}={}\\
  \shoveleft{{}=\PP_{k,1}{}\circ\mult_{1,k+1}\circ{}
    \Bigl((\xunderset{\mult_{1,k+1}}{1}
    - \xunderset{\mult_{k+1,1}^*}{3})(C^{k}\tensor C^{k+1})}\\*
  {}+ \xunderset{(-1)^{k}(\id\tensor\uu_{k,k}^{k-1,k})
    \circ\PP_{1,k+1}(C^{k}\tensor C^{k+1})}{3}
  + \xunderset{(\id\tensor\uu_{k,k}^{k-1,k})
  (C^{k+1}\tensor C^{k})}{2}\Bigr)\\*
  -\PP_{k,1}\circ(\idd{k}{}\tensor\bb')
  \circ(\xunderset{\mult_{1,k+1}}{1}
  - \xunderset{\mult_{k+1,1}^*}{3})(C^{k}\tensor C^{k+1})\\*
  +\xunderset{\mult_{k.k+1}(\bb'C^k\tensor C^{k+1})}{3}
  +\xunderset{(-1)^k\PP_{k,1}\circ\mult_{1,2}C^{k+1}\tensor\bb'C^k}{2}
  \\
  -\xunderset{\PP_{k,1}\circ\mult_{1,k+1}(C^k\tensor\bb'C^{k+1})}{1}
  -\xunderset{(-1)^k\PP_{k,1}\circ
    \mult_{1,k+1}(\bb'C^{k+1}\tensor C^k)}{2}=0,
\end{multline*}
where each of the three groups of terms (labeled with $1$, $2$,
and~$3$) vanish separately.  Equation~\eqref{work-91} is thus proved.

\subsubsection{} The subsequent calculations are straightforward but
quite tiresome.  To keep the presentation reasonably compact, we give
the details only for $\ell=1$.  This representatively illustrates the
general case because the quadratic term $\uu\circ\uu$ is already
present in the corresponding equation~\eqref{i=0} (compared
with~\eqref{ell-ge-2}, Eqs.~\eqref{i=0} involve an additional
complication due to the graded symmetrization
$F(C^k\tensor\bb'C^{k+1})+(-1)^{k}F(\bb'C^{k+1}\tensor C^k)$).

For $\ell=1$, the two equations that we must prove are given by
\begin{multline}\label{next-2}
  (\idd{(k-2)}\tensor\bb')
  \circ\uu_{k,k+1}^{k-2,k+2}(C^{k}\tensor C^{k+1})\\*
  +(-1)^k(\bb'\tensor\idd{(k+1)})\circ
  \uu_{k,k+1}^{k-1,k+1}(C^k\tensor C^{k+1})
  -\uu_{k-1,k+1}^{k-2,k+1}(\bb'C^k\tensor C^{k+1})\\*
  +(-1)^k\uu_{k,k}^{k-2,k+1}(C^k\tensor\bb'C^{k+1})  
  +\uu_{k,k}^{k-2,k+1}(\bb'C^{k+1}\tensor C^k)=\\
  {}=\uu_{k-1,k+1}^{k-2,k+1}\circ
  \uu_{k,k+1}^{k-1,k+1}(C^k\tensor C^{k+1})
\end{multline}
and
\begin{multline}\label{next-1}
  (-1)^k(\idd{(k-2)}\tensor\bb')
  \circ\uu_{k-1,k+2}^{k-2,k+2}(C^{k-1}\tensor C^{k+2})\\*
  {}+\uu_{k-1,k+1}^{k-2,k+1}(C^{k-1}\tensor\bb'C^{k+2})=0.
\end{multline}

To begin with~\eqref{next-1}, we recall recursive
relation~\eqref{recursion-one}.  We also use the recursive relation
\begin{equation}\label{it-also}
  \uu_{k-1,k+1}^{k-2,k+1}
  = (-1)^{k+1}\PP_{k-1,1}\circ\Bigl(\mult_{2,1}^* - \mult_{1,2}
  + (\id\tensor\uu_{k-1,k}^{k-2,k})
  \Bigr)\circ\PP_{1,k}
\end{equation}
that follows from~\thetag{$\boldsymbol{2k\!-\!2;k\!-\!2}$} (i.e., from
Eq.~\thetag{$\boldsymbol{2k;m}$} where we replace $k\to k-1$ and set
$m=k-2$).  We insert Eqs.~\eqref{recursion-one} and~\eqref{it-also}
in~\eqref{next-1} but keep $\uu_{k-1,k+1}^{k-2,k+1}$ ``unevaluated''
in the combination $\id\tensor\uu_{k-1,k+1}^{k-2,k+1}$ that arises
from~\eqref{recursion-one}.  As before, we use
identity~\eqref{obvious} (with $k\to k-1$).  Equation~\eqref{next-1}
that we must verify then becomes
\begin{multline}\label{work-92}
  (-1)^k\mult_{1,k}\circ(\mult_{2,1}^* - \mult_{1,2})\circ\PP_{1,k}
  (C^{k-1}\tensor C^{k+2})\\*
  +(-1)^{k+1}(\idd{(k-1)}\tensor\bb')
  \circ(\mult_{2,1}^* - \mult_{1,2})
  \circ\PP_{1,k}(C^{k-1}\tensor C^{k+2})\\*
  +(\idd{(k-1)}\tensor\bb')
  \circ(\id\tensor\uu_{k-1,k+1}^{k-2,k+1})
  \circ\PP_{1,k}(C^{k-1}\tensor C^{k+2})\\*
  -\mult_{1,k}\circ(\id\tensor\uu_{k-1,k+1}^{k-2,k+1})
  \circ\PP_{1,k}(C^{k-1}\tensor C^{k+2})\\*
  +(-1)^{k+1}(\mult_{2,1}^* - \mult_{1,2})
  \circ\PP_{1,k}(C^{k-1}\tensor\bb'C^{k+2})\\*
  +(-1)^{k+1}(\id\tensor\uu_{k-1,k}^{k-2,k})
  \circ\PP_{1,k}(C^{k-1}\tensor\bb'C^{k+2})=0.
\end{multline}

We now apply the induction hypothesis: for the \textit{preceding}
mapping $\uu_{k-1,k+1}^{k-2,k+1}$, Eqs.~\eqref{rest-eqs} imply that
\begin{equation}\label{hypo-1}
  (\idd{(k-2)}\tensor\bb')\circ\uu_{k-1,k+1}^{k-2,k+1}
  =-(-1)^k\uu_{k-1,k}^{k-2,k}\circ(\idd{(k-1)}\tensor\bb').
\end{equation}
Using this together with a simple identity
\begin{multline*}
  (\idd{(k-1)}\tensor\bb')\circ(\mult_{2,1}^* - \mult_{1,2})
  +(\mult_{1,k}-\mult_{k,1}^*)\circ(\idd{(k-1)}\tensor\bb')
  \circ\PP_{k,1}\\*
  -\mult_{1,k}\circ(\mult_{2,1}^* - \mult_{1,2})
  = \mult_{1,k}\circ\mult_{1,2}
  - \mult_{1,2}\circ\mult_{1,k+1},
\end{multline*}
we reduce Eq.~\eqref{work-92} to
\begin{multline}\label{work-95}
  (-1)^{k+1}(\mult_{1,k}\circ\mult_{1,2} - \mult_{1,2}\circ\mult_{1,k+1})
  -\mult_{1,k}\circ(\id\tensor\uu_{k-1,k+1}^{k-2,k+1})\\*
  -(-1)^k(\id\tensor\uu_{k-1,k}^{k-2,k})\circ\mult_{1,k+1}
  =0.
\end{multline}
Writing $\uu_{k-1,k+1}^{k-2,k+1} = \PP_{k-1,1}\circ\widetilde\uu$, we
have $\mult_{1,k}\circ(\id\tensor\uu_{k-1,k+1}^{k-2,k+1})=
\mult_{1,2}\circ(\id\tensor\widetilde\uu)$.  Recalling the actual form
of $\widetilde\uu$ from~\eqref{it-also}, we then see that
Eq.~\eqref{work-95} is identically satisfied (independently of any
properties of $\uu_{k-1,k}^{k-2,k}$).  \hbox{Hence, Eq.~\eqref{next-1}
  is proved}.

We next prove Eq.~\eqref{next-2}.  In its left-hand side,
$\uu_{k,k+1}^{k-2,k+2}$ is expressed from recursive
relation~\eqref{rec-second}.  We also express $\uu_{k,k+1}^{k-1,k+1}$
using Eq.~\eqref{rec-next}, except for the occurrences
of~$\uu_{k,k+1}^{k-1,k+1}$ originating from~\eqref{rec-second};
instead, we use Eq.~\eqref{ell=0} to rewrite the thus arising
combination $(\idd{(k-1)}\tensor\bb')\circ \uu_{k,k+1}^{k-1,k+1}$ as
\begin{multline*}
  (\idd{(k-1)}\tensor\bb')\circ
  \uu_{k,k+1}^{k-1,k+1}(C^{k}\tensor C^{k+1})={}\\
  =(-1)^{k}\uu_{k,k}^{k-1,k}(C^{k}\tensor\bb' C^{k+1})
  + \uu_{k,k}^{k-1,k}
  (\bb'C^{k+1}\tensor C^{k}).
\end{multline*}
In the right-hand side of~\eqref{next-2}, we use recursive
relations~\eqref{rec-next} and~\eqref{it-also} and then apply the
induction hypothesis, which implies the equation
\begin{multline}\label{prev-q}
  (-1)^k\uu_{k-1,k}^{k-2,k}
  (\bb'C^k\tensor C^k)
  +\uu_{k-1,k}^{k-2,k}\circ\uu_{k,k}^{k-1,k}(C^k\tensor C^k)\\*
  {}+(-1)^k(\idd{(k-2)}\tensor\bb'_{(k+1)})\circ\uu_{k,k}^{k-2,k+1}
  (C^k\tensor C^k)\\*
  + (\bb'_{(k-1)}\tensor\idd{k})\circ\uu_{k,k}^{k-1,k}(C^k\tensor C^k)
  =0,  
\end{multline}
which we now use to eliminate the term bilinear in~$\uu$ (as noted
above, equations evaluated on $C^k\tensor C^k$ require graded
symmetrization when the argument is stripped off).  Straightforward
manipulations involving~\eqref{obvious} and another obvious
identity,\enlargethispage{20pt}
\begin{equation*}
  (\bb'_{(k-1)}\tensor\idd{(k+1)})\circ\PP_{k,1}
  =\PP_{k-1,1}\circ(\id\tensor\bb'_{(k-1)}\tensor\idd{k})
\end{equation*}
then show that the terms arising from the second and third lines
in~\eqref{prev-q} cancel against similar terms in the left-hand side
of~\eqref{next-2}, and Eq.~\eqref{next-2} therefore becomes (where in
the hope to keep the derivation traceable, we have not yet used any
relations for $\uu_{k - 1, k + 1}^{ k - 2, k + 1}$ and
$\uu_{k,k}^{k-2,k+1}$)
\begin{multline}\label{mult-next}
  \uu_{k, k}^{k - 1, k}\bigl(
  C^k\tensor\bb'C^{k+1} + (-1)^{k} \bb'C^{k+1}\tensor C^{k}\bigr)\\*
  + \mult_{1, k}\circ\Bigl(
  (\id\tensor\uu_{ k - 1, k + 1}^{ k - 2, k + 1})
  (C^k\tensor C^{k+1}) - 
  (\id\tensor\uu_{ k, k}^{ k - 2, k + 1})
  (C^{k+1}\tensor C^k)\\*
  \shoveright{{}- (-1)^k (\id\tensor\uu_{ k, k}^{ k - 2, k + 1})
  \circ\PP_{1,k+1} (C^k\tensor C^{k+1})\Bigr)}\\*
  -   (\idd{(k - 1)}\tensor\bb')\circ
  (\id\tensor\uu_{ k - 1, k + 1}^{ k - 2, k + 1})
  (C^k\tensor C^{k+1})
  - \PP_{1,k-1}\circ\uu_{k - 1, k + 1}^{ k - 2, k + 1}
  (\bb'C^k\tensor C^{k+1})\\*
  + (-1)^k (\id\tensor\bb'_{(k - 1)}\tensor\idd{k})
  \circ(\mult_{1, k + 1}-\mult_{k + 1, 1}^*) (C^k\tensor C^{k+1})\\*
  \shoveright{
    {}+ (-1)^k \PP_{1,k-1}\circ\uu_{k, k}^{k - 2, k + 1}
    (C^k\tensor\bb'C^{k+1})
    +\PP_{1,k-1}
    \circ\uu_{k, k}^{ k - 2, k + 1} (\bb'C^{k+1}\tensor C^{k}) ={}}\\
  = (-1)^{k}
  \mult_{1,2}\circ\Bigl(
  (\mult_{1, k + 1} - \mult_{k + 1, 1}^*)(C^k\tensor C^{k+1})
  + (\id\tensor\uu_{ k, k}^{ k - 1, k}) (C^{k+1}\tensor C^k)\\*
  \shoveright{{}+ (-1)^k (\id\tensor\uu_{ k, k}^{ k - 1, k})\circ
  \PP_{1,k+1} (C^k\tensor C^{k+1})\Bigr)}\\*
   + (-1)^{k + 1}(\id\tensor\uu_{ k - 1, k}^{ k - 2, k})\circ
  (\mult_{1, k + 1} - \mult_{k + 1, 1}^*)(C^k\tensor C^{k+1}) \\*
  + 2\bigl(\id\tensor(\uu_{ k - 1, k}^{ k - 2, k}\circ
  (\bb'_{(k)}\tensor\idd{k}))\bigr)\circ
  \Sym{k}{2, k + 1;k + 2, 2k + 1}
  (C^{k+1}\tensor C^k).
\end{multline}
The $\Sym{k}{2, k + 1; k + 2, 2k + 1}$ symmetrization arises here
because of the symmetry of the argument in~\eqref{prev-q}.

Next, using the recursive relations for $\uu_{k - 1, k + 1}^{ k - 2, k
  + 1}$ and $\uu_{k,k}^{k-2,k+1}$ in the left-hand side
of~\eqref{mult-next}, we obtain in accordance with~\eqref{it-also}
that
\begin{multline}\label{work-100}
  {-}\PP_{1,k-1}\circ\uu_{k - 1, k + 1}^{ k - 2, k + 1}
  (\bb'C^k\tensor C^{k+1})=\\
  {}=(-1)^k(\mult_{1,k}-\mult_{k,1}^*)(\bb'C^k\tensor C^{k+1})
  {}+(-1)^k(\id\tensor\uu_{k-1,k}^{k-2,k})
  \circ\PP_{1,k}(\bb'C^k\tensor C^{k+1}),
\end{multline}
where the last term is readily seen to cancel one of the two terms
produced by the $\Sym{k}{2, k + 1;k + 2, 2k + 1}$ symmetrization in
the right-hand side of~\eqref{mult-next}.  Next, recalling the
recursive relation
\begin{multline}\label{rec-final}
  \uu_{k,k}^{k-2,k+1}(C^k\tensor C^k)=\\
  {}=-\PP_{k-1,1}\circ(\id\tensor\uu_{k-1,k}^{k-2,k})(C^k\tensor C^k)
  -(-1)^k\PP_{k-1,1}\circ\uu_{k,k}^{k-1,k}(C^k\tensor C^k)
\end{multline}
(which is Eq.~\eqref{uk1k1metc} for $k\to k-1$ and $m=k-2$), we see
that the terms in the right-hand side cancel other occurrences of
$\uu_{k-1,k}^{k-2,k}$ and~$\uu_{k,k}^{k-1,k}$ in~\eqref{mult-next}.
In the left-hand side of~\eqref{mult-next}, we also evaluate
\begin{multline}
  (-1)^k (\id\tensor\bb'_{(k - 1)}\tensor\idd{k})
  \circ(\mult_{1, k + 1}-\mult_{k + 1, 1}^*) (C^k\tensor C^{k+1})\\
  =(-1)^k(\mult_{1,k}-\mult_{k,1}^*)
  \bigl(((\id\tensor\bb')C^k)\tensor C^{k+1}\bigr),
\end{multline}
which cancels most of the first term in the right-hand side
of~\eqref{work-100}.  For the combination $(\idd{(k -
  1)}\tensor\bb')\circ (\id\tensor\uu_{ k - 1, k + 1}^{ k - 2, k +
  1})$ in~\eqref{mult-next}, we again use Eq.~\eqref{hypo-1} implied
by the induction hypothesis.

Equation~\eqref{mult-next} thus becomes
\begin{multline}\label{mult-next2}
  \mult_{1, k}\circ\Bigl(
  (\id\tensor\uu_{ k - 1, k + 1}^{ k - 2, k + 1})
  (C^k\tensor C^{k+1}) - 
  (\id\tensor\uu_{ k, k}^{ k - 2, k + 1})
  (C^{k+1}\tensor C^k)\\*
  {}- (-1)^k (\id\tensor\uu_{ k, k}^{ k - 2, k + 1})
  \circ\PP_{1,k+1} (C^k\tensor C^{k+1})\Bigr)\\*
  \shoveright{{} +
    \underline{(-1)^k(\mult_{1,k}-\mult_{k,1}^*)
      (\mult_{1,2}C^k\tensor C^{k+1})}={}}\\
  {}= (-1)^{k}
  \mult_{1,2}\circ\Bigl(
  \underline{(\mult_{1, k + 1} - \mult_{k + 1, 1}^*)(C^k\tensor C^{k+1})}
  + (\id\tensor\uu_{ k, k}^{ k - 1, k}) (C^{k+1}\tensor C^k)\\*
  {{}    
    + (-1)^k (\id\tensor\uu_{ k, k}^{ k - 1, k})\circ
    \PP_{1,k+1} (C^k\tensor C^{k+1})\Bigr)}\\*
   + (-1)^{k + 1}(\id\tensor\uu_{ k - 1, k}^{ k - 2, k})\circ
  (\mult_{1, k + 1} - \mult_{k + 1, 1}^*)(C^k\tensor C^{k+1}) \\*
  +(\id\tensor\uu_{k-1,k}^{k-2,k})\circ\mult_{1,2}
  (C^{k+1}\tensor C^k).
\end{multline}
It can be immediately seen that the underlined terms can be replaced
with
\begin{equation}\label{m-remain}
  (-1)^k(\mult_{1,k}\circ\mult_{1,2} - \mult_{1,2}\circ\mult_{1,k+1})
\end{equation}
in the left-hand side.  Next, recursive relation~\eqref{it-also}
allows us to calculate
\begin{multline}\label{elegant}
  \mult_{1,k}\circ(\id\tensor\uu_{k-1,k+1}^{k-2,k+1})=\\
  \begin{aligned}
    &=(-1)^{k+1}\mult_{1,k}\circ
    \bigl(\id\tensor\bigl(\PP_{k-1,1}\circ(\mult_{2,1}^* - \mult_{1,2}
    + (\id\tensor\uu_{k-1,k}^{k-2,k}))\circ\PP_{1,k}\bigr)\bigr)\\
    &=(-1)^{k+1}
    \mult_{1,k}\circ
    \PP_{k,2}\circ\Bigl(\mult_{3,2}^* - \mult_{2,3}
    + (\idd{2}\tensor\uu_{k-1,k}^{k-2,k})
    \Bigr)\circ\PP_{2,k+1}\\
    &=(-1)^{k+1}(\mult_{1,k}\circ\mult_{1,2} -
    \mult_{1,2}\circ\mult_{1,k+1})
    +(-1)^{k+1}(\id\tensor\uu_{k-1,k}^{k-2,k})\circ\mult_{1,k+1}.
  \end{aligned}
\end{multline}
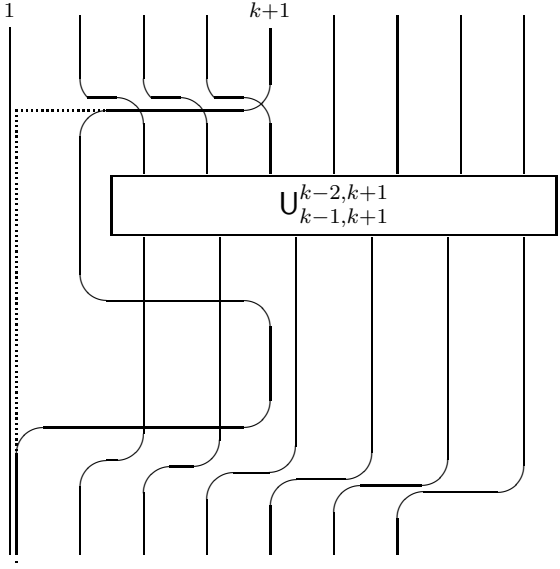
\begin{SCfigure}[60][tb]
  \begin{xy}
    <2pc,0pc>:(4,8.5)*{\mbox{${}^{k+1}$}}
    \PATH~={**\dir{-}?>}
    `d (0,7)
    `  (1,6)
    `  (4,4)
    `  (4,2)
    `  (0,2)
    `  (0,0)
    '  (0,0)
    \POS(3,7)
    \PATH~={**\dir{.}?>}
      ' (0,7)
      ' (0,0)
    \POS(-.1,8.5)*{\mbox{${}^{1}$}}
    \PATH~={**\dir{-}?>}
    ' (-.1,0)
    \POS(1,8.5)
    \PATH~={**\dir{-}?>}
    ' (1,7.5)
    `^dr (2,7.5)
    `r (2,6)
    '(2,6)
    \POS(2,8.5)
    \PATH~={**\dir{-}?>}
    ' (2,7.5)
    `^dr (3,7.5)
    `r (3,6)
    '(3,6)
    \POS(3,8.5)
    \PATH~={**\dir{-}?>}
    ' (3,7.5)
    `^dr (4,7.5)
    `r (4,6)
    '(4,6)
    \POS(5,8.5)
    \PATH~={**\dir{-}?>}
    ' (5,6)
    \POS(6,8.5)
    \PATH~={**\dir{-}?>}
    ' (6,6)
    \POS(7,8.5)
    \PATH~={**\dir{-}?>}
    ' (7,6)
    \POS(8,8.5)
    \PATH~={**\dir{-}?>}
    ' (8,6)
    \POS(5,5.5)
    *+\txt{\framebox{\kern60pt$\uu_{k-1,k+1}^{k-2,k+1}$\kern60pt}}="uu"
    \POS(2,5)
    \PATH~={**\dir{-}?>}
    ' (2,1.9)
    ` (1,1.9)
    ` (1,1)
    ' (1,0)
    \POS(3.2,5)
    \PATH~={**\dir{-}?>}
    ' (3.2,1.8)
    ` (2,1.8)
    ` (2,0)
    ' (2,0)
    \POS(4.4,5)
    \PATH~={**\dir{-}?>}
    ' (4.4,1.7)
    ` (3,1.7)
    ` (3,0)
    ' (3,0)
    \POS(5.6,5)
    \PATH~={**\dir{-}?>}
    ' (5.6,1.6)
    ` (4,1.6)
    ` (4,0)
    ' (4,0)
    \POS(6.8,5)
    \PATH~={**\dir{-}?>}
    ' (6.8,1.5)
    ` (5,1.5)
    ` (5,0)
    ' (5,0)
    \POS(8.0,5)
    \PATH~={**\dir{-}?>}
    ' (8.0,1.4)
    ` (6,1.4)
    ` (6,0)
    ' (6,0)
  \end{xy}\kern-40pt
  \caption[Graphical representation of braiding and multiplication
  operators]{\small Graphical representation of
    $\protect\mult_{1,k}\circ
    \protect\PP_{k,2}\circ(\idd{2}\tensor\uu_{k-1,k}^{k-2,k})
    \circ\protect\PP_{2,k+1}$ (depicted for $k=4$) showing that this
    expression coincides with
    $(\id\tensor\uu_{k-1,k}^{k-2,k})\circ\protect\mult_{1,k+1}$ (the
    dotted line representation for the $(k+1)$th tensor factor).  The
    operators $\protect\PP_{k,2}$ and $\protect\mult_{1,k}$ are given
    by braiding and joining the strings representing the corresponding
    tensor factors.  The order of applying operations from right to
    left here corresponds to the order from top down.}
  \label{illustration}
\end{SCfigure}
\mbox{}\kern-2pt Again, these manipulations do not rely on any
properties of $\uu_{k-1,k}^{k-2,k}$, see Fig.~1 for the corresponding
rearrangement of the term involving $\uu_{k-1,k}^{k-2,k}$; the terms
involving $\mult_{3,2}^*$ and $\mult_{2,3}$ are transformed similarly.
The last term cancels a similar term in the right-hand side
of~\eqref{mult-next2}, the $\mult\circ\mult$ terms cancel those
in~\eqref{m-remain}, and Eq.~\eqref{mult-next2} therefore becomes
\begin{multline*}
  -\mult_{1, k}\circ\Bigl(
  (\id\tensor\uu_{ k, k}^{ k - 2, k + 1})
  (C^{k+1}\tensor C^k)
  + (-1)^k (\id\tensor\uu_{ k, k}^{ k - 2, k + 1})
  \circ\PP_{1,k+1} (C^k\tensor C^{k+1})\Bigr)\\*
  {}= 
  \mult_{1,2}\circ\Bigl(  
  (-1)^{k}(\id\tensor\uu_{ k, k}^{ k - 1, k}) (C^{k+1}\tensor C^k)
  + (\id\tensor\uu_{ k, k}^{ k - 1, k})\circ
  \PP_{1,k+1} (C^k\tensor C^{k+1})\Bigr)\\*
  + (-1)^{k}(\id\tensor\uu_{ k - 1, k}^{ k - 2, k})\circ
  \mult_{k + 1, 1}^*(C^k\tensor C^{k+1})
  +(\id\tensor\uu_{k-1,k}^{k-2,k})\circ\mult_{1,2}
  (C^{k+1}\tensor C^k).
\end{multline*}
This is readily seen to be an identity in view of recursive
relation~\eqref{rec-final}, via an argument entirely similar to the
one applied to~\eqref{elegant} above.  Equation~\eqref{next-2} is
proved.

\subsection{The end of the proof of Theorem~\ref{thm:Omega2}}
\label{sec:end-proof-thm}
\addcontentsline{toc}{section}{\numberline{\thesubsection.}The end of
  the proof of Theorem~\ref{thm:Omega2}} It remains to prove
Eqs.~\eqref{C3}, or $\Xterm^{(3)}=0$, for the BRST
differential~$\BRST$.  The proof involves algebraic consequences of
Eqs.~\eqref{C2}, which are already established.  Some details are as
follows.  We write $\BRST$ as in~\eqref{gen-reducible} in the case
where $\tau_A$ and $U_{AB}$ have the form (with summations over
repeated indices)
\begin{equation}\label{x-generators}
  \begin{split}    
    \tau_A={}& t_A+Z_A^B\ghP_B,\\
    U_{AB}={}&U_{AB}^C\ghP_C+\thalf(-1)^{\varepsilon_C}
    U_{AB}^{CD}\ghP_D\ghP_C
  \end{split}
\end{equation}
with $\ghP$- and $C$-independent $t_A$, $Z_A^B$, $U_{AB}^C$, and
$U_{AB}^{CD}$.  Here, $\Grass{C^A}=\Grass{\ghP_A}=\varepsilon_A+1$.
In general, the $\oZ_2$-gradings are given by
$\Grass{t_A}=\varepsilon_A$,
$\Grass{Z_A^B}=\varepsilon_A+\varepsilon_B+1$,
$\Grass{U_{AB}^C}=\varepsilon_A+\varepsilon_B+\varepsilon_C$, and
$\Grass{U_{AB}^{CD}}=\varepsilon_A+\varepsilon_B+\varepsilon_C
+\varepsilon_D+1$, but in our specific case, these are all even.
Expanding Eqs.~\eqref{C1}--\eqref{C3} in the ghost momenta, we then,
in addition to the basic relations $Z_A^B t_B = 0$, $Z_A^B Z_B^C=0$,
and $[t_A,t_B]=U_{AB}^C t_C$, obtain the commutators involving
$Z_A^B$,
\begin{multline}\label{C2P1}
  [t_A,Z_B^D]-(-1)^{\varepsilon_A\varepsilon_B}[t_B,Z_A^D]={}\\
  {}=U_{AB}^C Z_C^D - U_{AB}^{DC} t_C
  - (-1)^{\varepsilon_B}Z_A^C U_{CB}^D
  + (-1)^{\varepsilon_A(1+\varepsilon_B)}Z_B^C U_{CA}^D
  + \thalf U_{AB}^{EC} U_{CE}^D,
\end{multline}
\begin{multline}\label{C2P2}
  (-1)^{(\varepsilon_B+1)(\varepsilon_D+1)}[Z_A^D, Z_B^E]
  -(-1)^{(\varepsilon_B+1)(\varepsilon_E+1)
    +\varepsilon_D\varepsilon_E}[Z_A^E, Z_B^D]={}\\
  {}=(-1)^{\varepsilon_D}U_{AB}^{DC}Z_{C}^E
  -(-1)^{\varepsilon_E(1+\varepsilon_D)}U_{AB}^{EC}Z_C^D\\*
  +(-1)^{\varepsilon_B}Z_A^C U_{CB}^{DE}
  -(-1)^{\varepsilon_A(1+\varepsilon_B)}Z_B^C U_{CA}^{DE}
  -\thalf U_{AB}^{FC}U_{CF}^{DE},
\end{multline}
and (from~\eqref{C3}) the Jacobi identity
\begin{equation}\label{C3P1}
  (-1)^{\varepsilon_A\varepsilon_C}U_{AB}^DU_{DC}^E +
  \mathrm{cycle}(A,B,C)=0
\end{equation}
and its generalizations involving the four-indexed $U_{AB}^{CD}$,
\begin{multline}\label{C3P2}
  (-1)^{\varepsilon_A\varepsilon_C}
  \Bigl(U_{AB}^DU_{DC}^{EF}+
  (-1)^{(\varepsilon_E+1)\varepsilon_C}U_{AB}^{ED}U_{DC}^F\\* 
  - (-1)^{(\varepsilon_F+1)\varepsilon_C +
    \varepsilon_E\varepsilon_F}U_{AB}^{FD}U_{DC}^E\Bigr) 
  {}+\mathrm{cycle}(A,B,C)=0
\end{multline}
and
\begin{multline}\label{C3P3}
  (-1)^{(\varepsilon_E+1)(\varepsilon_C+1) + \varepsilon_A\varepsilon_C
    + \varepsilon_E\varepsilon_G}U_{AB}^{ED}U_{DC}^{FG}\\*
  +\mathrm{cycle}(A,B,C) + \mathrm{cycle}(E,F,G)=0.
\end{multline}

The tedious proof in~\bref{sec:proof}--\bref{sec:end-proof-eq} is
nothing but the demonstration of Eqs.~\eqref{C2P1} and~\eqref{C2P2}.
The remaining equations~\eqref{C3P1}--\eqref{C3P3} can be shown by
studying the algebraic consistency conditions for~\eqref{C2P1}
and~\eqref{C2P2} and using specific properties of~$Z_A^B$ and
``selection rules'' for $U_{AB}^C$ and~$U_{AB}^{CD}$.  For this, we
evaluate the left-hand side of the Jacobi identity
$(-1)^{\varepsilon_A\varepsilon_C}[[\tau_A,\tau_B],\tau_C]+
\mathrm{cycle}(A,B,C)=0$.  The double commutator involves terms of
zeroth, first, and second orders in the ghost momenta,
\begin{equation*}
  [[\tau_A,\tau_B],\tau_C]=J_{ABC} + J_{ABC}^H\ghP_H +
  J_{ABC}^{HG}\ghP_G\ghP_H,
\end{equation*}
where
\begin{equation*}
  J_{ABC}=U_{AB}^DU_{DC}^E t_E
\end{equation*}
and
\begin{multline*}
  J_{ABC}^H=-(-1)^{\varepsilon_C}U_{AB}^DZ_D^EU_{EC}^H
  +(-1)^{\varepsilon_D(\varepsilon_C+1)}U_{AB}^DZ_C^EU_{ED}^H
  + U_{AB}^DU_{DC}^EZ_E^H\\*
  + \thalf U_{AB}^DU_{DC}^{EK}U_{KE}^H
  +(-1)^{\varepsilon_E\varepsilon_H}U_{AB}^DU_{DC}^{EH}t_E
  +(-1)^{\varepsilon_D+(\varepsilon_H+1)(\varepsilon_D+\varepsilon_C)}
  U_{AB}^{DH}U_{DC}^Et_E\\*
  +(-1)^{\varepsilon_A\varepsilon_C
    + \varepsilon_B\varepsilon_C + \varepsilon_B}[t_C,Z_A^D]U_{DB}^H
  -(-1)^{\varepsilon_A\varepsilon_B + \varepsilon_B\varepsilon_C
    + \varepsilon_A\varepsilon_C + \varepsilon_A}[t_C,Z_B^D]U_{DA}^H
\end{multline*}
(the expression for $J_{ABC}^{HG}$ is totally straightforward but
rather unwieldy).  The vanishing of the cyclic sum of
$(-1)^{\varepsilon_A\varepsilon_C} J_{ABC}$ does not allow us to
conclude immediately that the cyclic sum of
$(-1)^{\varepsilon_A\varepsilon_C} U_{AB}^DU_{DC}^E$ vanishes, because
as follows from comparing~\eqref{Omega-A} and the first equation
in~\eqref{x-generators}, the $t_A$ term is present in $\tau_A$ only in
the ``irreducible sector,'' where $A=a$, in which case $t_a$ are
elements of a basis in~$\algA$.  We must therefore distinguish between
the cases where the indices $A,B,C$ take the ``lower'' and ``higher''
values, i.e., label elements of a basis in~$\algA$ and elements of a
basis in $\algA^{\tensor n}$ with $n\geq2$, respectively.  But when
all three indices take ``higher'' values, the commutator terms in the
formula for $J_{ABC}^H$ vanish.  Moreover, selection rules for the
$\uu$ mappings (that $\uu_{m,n}^{m',n'}$ can be nonzero only for
$m+n=m'+n'+1$) show that the terms involving $t_E$ can also be
discarded for generic ``higher'' values of $A,B,C$.  Among the
remaining four terms in $J_{ABC}^H$, one is proportional to the unit
element of the algebra, while the other three are not (when evaluated
on generic elements in the corresponding tensor powers of~$\algA$).
To these three elements (in the first line in the last formula), we
apply the same trick of taking the quotient over $\oC\cdot1$ as
in~\bref{sec:recursion}.  In the current normalization, this amounts
to replacing $Z_A^B$ with $(-1)^{\varepsilon_A}\delta_A^B$.  Two of
the three terms then cancel each other, while the cyclic sum of the
remaining one gives the desired relation~\eqref{C3P1};
% \begin{equation*}
%   (-1)^{\varepsilon_A\varepsilon_C} U_{AB}^D U_{DC}^H
%   + \mathrm{cycle}(A,B,C)=0.
% \end{equation*}
% This relation is a part of Eq.~\eqref{C3} (and is obtained by
% expanding the latter in the ghost momenta); 
% the remaining part of~\eqref{C3} is shown similarly.
the remaining equations are shown similarly.

\section{Additional Remarks}\label{sec:other}
\subsection{``Hamiltonians,''  cohomology, and observables} 
\addcontentsline{toc}{section}{\numberline{\thesubsection.}Hamiltonians,
  cohomology, and observables} 
The cohomology problem for the differential $\BRST$ can in principle
be considered on different spaces.  As with the general form
of~$\BRST$, we borrow one such setting from quantization of gauge
theories.  In genuine gauge theories, the (quantum) Hamiltonian
$\cH\in\AQgh$ commutes with $\BRST$ (and is therefore naturally
considered modulo $\BRST$-exact terms).  A physical requirement is
$\gh\cH=0$.  In closed-algebra gauge theories, Hamiltonians are
typically at most linear in the~$C$ ghosts but can depend on the ghost
momenta.  We restrict ourself to the at-most-bilinear dependence on
the ghost momenta and moreover choose the part depending on the ghost
momenta to be ``scalar'' in the sense that it commutes with elements
of~$\AQ$ (the ``original'' operators of the theory, not involving
ghosts).  In our setting, such Hamiltonians are therefore given~by
\begin{equation*}
  \cH=H+\sum_{r\geq1}\Contr{\vv_r(C^r)}{P_r}
  +\mathop{\sum\!\sum}_{s>r\geq1}
  \Contr{\vv_{r+s}^{r,s}(C^{r+s})}{P_r\tensor P_s},
\end{equation*}
where $H\in\algA$ and
\begin{align*}
  \vv_r\in{}&\Hom(\algA^{\tensor r},\algA^{\tensor r}),\\
  \vv_m^{r,s}\in{}&\Hom(\algA^{\tensor m},
  \algA^{\tensor r}\tensor\algA^{\tensor s}).
\end{align*}
The terms that are linear in the~$C$ ghosts in the commutator are
given by
\begin{multline*}
  [\cH,\BRST]\Bigr|_{C^{\tensor1}}
  =[H,C^1] +\vv_1(C^1)\\*
  + \sum_{s\geq1}\contr{\ad{1}{H}\ldot\Z{s+1}(C^{s+1})}{P_s}
  + \sum_{s\geq2}(-1)^s \contr{\vv_{s+1}^{1,s}(C^{s+1})}{P_s}\\*
  + \sum_{s\geq1}(-1)^s\contr{\Z{s+1}\circ\vv_{s+1}(C^{s+1})}{P_s}
  + \sum_{s\geq1}(-1)^s\contr{(\id\tensor\vv_s)\circ\Z{s+1}(C^{s+1})}{
    P_s}\\*
  +\mathop{\sum\!\sum}_{s>r\geq1}(-1)^{rs+r+s}
  \Contr{\uu_{r,s}^{r+s-1}\circ\vv_{r+s}^{r,s}(C^{r+s})}{P_{r+s-1}}\\*
  + \sum_{r\geq1}\!\sum_{s\geq3}(-1)^{r(s+1)}
  \contr{
    (\Z{r+1}\tensor\idd{s})\circ\vv_{r+s+1}^{r+1,s}(C^{r+s+1})}
  {P_r\tensor P_s}\\*
  + \sum_{r\geq1}\!\sum_{s\geq1}
  (-1)^s
  \contr{\PP_{1,r+1}(\idd{r}\tensor\Z{s+1})\circ
   \vv_{r+s+1}^{r,s+1}(C^{r+s+1})}{P_r\tensor P_s}\\*
  + \mathop{\sum\!\sum}_{s>r\geq1}\!\!
  \sum_{m=1}^{\intpart{\frac{r+s-2}{2}}}
  (-1)^{rs+r+s}
  \Contr{\uu_{r,s}^{m,r+s-m-1}\circ\vv_{r+s}^{r,s}(C^{r+s})}{
    P_m\tensor P_{r+s-m-1}}.
\end{multline*}

For any $H\in\algA$, there exist mappings $\vv_r$ 
such that $[\cH,\BRST]=0$.\footnote{We note that from the standpoint
  of constrained systems, a ``Hamiltonian'' $H\in\algA$ is a linear
  combination of constraints (elements of any chosen basis
  in~$\algA$).  The associative-algebra counterparts of more general
  Hamiltonians can be realized by taking $H$ to be an arbitrary
  endomorphism of a representation space of~$\algA$.}  Indeed, the
vanishing of $P_n$-independent terms in the commutator is expressed~as
\begin{equation*}
  [H,C^1] +\vv_1(C^1)=0,
\end{equation*}
and the vanishing of the terms linear in~$P_n$ as
\begin{equation*}  
  \ad{1}{H}\ldot\Z{n+1}
  + (-1)^n\Z{n+1}\circ\vv_{n+1}
  + (-1)^n(\id\tensor\vv_n)\circ\Z{n+1}=0,
\end{equation*}
which immediately yields a solution
\begin{equation*}
  \vv_{n}=(-1)^n\Lie_{H},
\end{equation*}
and it remains to note that $\Lie_H\circ\uu_{m,n}^{\ell,m+n-1-\ell}
-\uu_{m,n}^{\ell,m+n-1-\ell}\circ\Lie_H=0$ for any~$H\in\algA$.
Moreover, $\cH$ of the above form is not $\BRST$-exact, and the
cohomology of~$\BRST$ on the space of such operators therefore
contains the algebra~$\algA$.  The cohomology of $\BRST$ on objects
with other, nonzero ghost numbers determines what is called
observables in the standard quantum setting.

\subsection{Automorphisms} 
\addcontentsline{toc}{section}{\numberline{\thesubsection.}Automorphisms}
Obviously, the relation $\BRST\,\BRST=0$ is preserved by similarity
transformations $\BRST\mapsto \Utransf\,\BRST\,\Utransf^{-1}$, where
$\Utransf$ is an arbitrary invertible operator on the bar
resolution; %\pagebreak[3] 
a natural subclass of such operators is given by ``inner''
automorphisms with $\Utransf=e^{\Gtransf}$ (assuming that the
exponential mapping exists), where $\Gtransf$, with $\gh \Gtransf=0$,
is an arbitrary operator on the bar resolution.  Following the pattern
set above, we can further restrict the class of transformations by
requiring $\Gtransf$ to have scalar coefficients except in the
``constant'' term (even with these operators, the transformed BRST
differential would no longer manifestly look like the one
corresponding to a closed algebra).  In particular, this class
contains transformations described in~\cite{[BT-UUU]} that act
``covariantly'' on the multiplication
$\mult\in\Hom(\algA\tensor\algA,\algA)$, thereby ``preserving
associativity.''  Restricting to transformations that preserve the
maximum powers in the ghost and ghost momenta expansion, it may be
possible to obtain solutions for the $\uu$ mappings in which the
mappings listed in~\eqref{eq:vanishing} do not necessarily vanish.

\subsection{Weyl ordering}
\addcontentsline{toc}{section}{\numberline{\thesubsection.}Weyl
  ordering} The $CP$ ordering of the ghost operators (ghosts to the
left and momenta to the right) was chosen in constructing the BRST
differential, but other orderings are often preferred in ``genuine''
quantum systems; most popular are the Weyl (totally symmetric) and
Wick orderings.  In application to associative algebras, the $CP$
ordering is ``nonminimal'' in the following sense.  As noted above, a
totally straightforward part of the construction of~$\BRST$ is related
to the (graded) \textit{antisymmetrized} part of the multiplication
in~$\algA$, i.e., to the associated Lie algebra~\eqref{fabc-skew}.
The entire reducible gauge theory formalism then serves to incorporate
the (graded) \textit{symmetrized} part of the multiplication,
$\{t_a,t_b\}=\sum_c\,f_{\{ab\}}^c\,t_c$.  But the antisymmetric sector
is actually ``admixed'' to the relations in the ``reducible'' part of
the formalism; for example, the tensor $\Z{}_{ab}^c$
(see~\eqref{Z-zero}) involves an antisymmetric part with respect to
its lower indices, and the corresponding projection of the equation
for $\Z{2}$ is satisfied as a consequence of~\eqref{fabc-skew}.

In the Weyl ordering, this redundancy is eliminated, but most of the
relations following from $\BRST\BRST=0$ become somewhat less
transparent.  Using $Z^{c}_{ab}$, etc.\ to denote the coefficients in
the Weyl-ordered BRST operator, we then have
\begin{equation*}
  \sum_{c}(Z^c_{ab}\,t_c + t_c\,Z^c_{ab})=0
\end{equation*}
instead of Eq.~\eqref{Z-zero}, and hence,
\begin{equation*}
  Z^c_{ab}=t_a\,\delta_b^c-\thalf f^c_{\{ab\}}.
\end{equation*}
Next, instead of~\eqref{higher-Z-zero} with $n=2$, we have an equation
with a nonzero right-hand side,
\begin{equation*}
  Z_{def}^{ab}\,Z_{ab}^c + Z_{ab}^c\,Z_{def}^{ab}
  =-\thalf f_{[df]}^b\,f_{[eb]}^c
\end{equation*}
with
\begin{equation*}
  Z_{def}^{ab}=t_d\,\delta_e^a\delta_f^b -
  \thalf\,f_{\{de\}}^a\delta_f^b
  +\thalf\,\delta_d^a f_{\{ef\}}^b,
\end{equation*}
and so on.

We also note that in the notation similar to that
in~\bref{sec:closed}, with $\tau_A$ and $U_{AB}$ now used to denote
the respective coefficients in the Weyl form of the BRST
operator~\eqref{gen-reducible} for a closed-algebra theory, the
equations analogous to~\eqref{C1}--\eqref{C3} (derived under the same
conditions as in~\eqref{scalar-case}) begin with \textit{zero}-order
ones in the expansion in the ghosts (with summations over repeated
indices understood),
\begin{equation*}
  (-1)^{\Grass{C^A}+1}[[\tau_A,C^B],[\tau_B,C^A]]
  +(-1)^{\Grass{C^B} + \Grass{C^D}}
  [[\tau_D,C^B],C^A]\,[U_{AB},C^D] = 0.
\end{equation*}
Equations corresponding to the first and second order in the ghosts
are given by
\begin{multline*}
  [\tau_D,C^B]\tau_B
  +(-1)^{(\Grass{C^B}+1)(\Grass{C^D}+1)}\tau_B[\tau_D,C^B]\\*
  -(-1)^{\Grass{C^B}}\tfrac{1}{4}[[U_{DE},C^B],C^A][U_{AB},C^E]=0
\end{multline*}
and
\begin{multline*}
  [\tau_A,\tau_B]=[U_{AB},C^D]\tau_D
  + (-1)^{\Grass{C^B}}[\tau_A,C^D]U_{DB}\\*
  + (-1)^{(\Grass{C^A}+1)\Grass{C^B}}[\tau_B,C^D]U_{DA}.
\end{multline*}
In a theory with a closed algebra, the $C^3$-equations are not
modified compared with~\eqref{C3}.

\subsection*{Acknowledgments} We thank B.L.~Feigin, M.A.~Grigoriev, and
I.V.~Tyutin for useful discussions.  This paper was supported in part
by the RFBR (grants 02-01-00930 and 00-15-96566), grant
LSS-1578.2003.2, INTAS (grant 00-00262), and the Foundation for
Promotion of Russian Science.

\end{document}